\documentclass[a4paper]{article}
\usepackage{amsfonts}
\usepackage{amssymb}
\newtheorem{st}   {Theorem}   [section] 
\newtheorem{prop}[st]{Proposition}
\newtheorem{lem}[st] {Lemma}   
\newtheorem{cor}[st] {Corollary} 
    
\newtheorem{defin}[st]    {Definition}
\newtheorem{exam}[st]{Example} 
\newtheorem{num}[st]{}

\newcommand{\lie}{\ensuremath{\frak{g}}}
\newcommand{\B}{\ensuremath{\mathcal{B}}}
\newcommand{\A}{\ensuremath{\mathcal{A}}}
\newcommand{\F}{\ensuremath{\mathcal{F}}}

\newcommand{\Ha}{\ensuremath{\mathcal{H}}}
\newcommand{\xx}{\ensuremath{\mathcal{X}}}

\newcommand{\rmap}{\longrightarrow}
\newcommand{\Boxe}{\raisebox{.8ex}{\framebox}}

\newcommand{\del}{\partial}
\newcommand{\ps}{{\raise 1pt\hbox{\tiny (}}}

\newcommand{\pss}{{\raise 1pt\hbox{\tiny [}}}
\newcommand{\pdd}{{\raise 1pt\hbox{\tiny ]}}}
\newcommand{\pd}{{\raise 1pt\hbox{\tiny )}}}

\newcommand{\bs}{{\raise 1pt\hbox{\tiny [}}}
\newcommand{\bd}{{\raise 1pt\hbox{\tiny ]}}}
\def\cross{\mathinner{\mathrel{\raise0.8pt\hbox{$\scriptstyle>$}}
                 \joinrel\mathrel\triangleleft}}

\def\compose{{\raise 1pt\hbox{$\scriptscriptstyle\circ$}}}
\def\dcross{{\raise 0.5pt\hbox{$\scriptscriptstyle\boxtime$}}}

\usepackage[all]{xy}
\usepackage{epsf} 
\usepackage{amsfonts}
\usepackage{amssymb}

\begin{document}

\title{Cyclic cohomology of Hopf algebras, and a non-commutative Chern-Weil
  theory
\thanks{Research supported by NWO}}
\author { by Marius Crainic}
\date{Utrecht University, Department of Mathematics,\\
P.O.Box:80.010,3508 TA Utrecht, The Netherlands,\\ 
e-mail: crainic@math.uu.nl
}
\pagestyle{myheadings}
\maketitle
\begin{abstract} We give a motivation/interpertation of Connes-Moscovici's cyclic cohomology of Hopf algebras in the spirit of the Cuntz-Quillen formalism.  Furthermore, we introduce a non-commutative Weil complex, which connects the work of Gelfand and Smirnov with cyclic cohomology. We show that the Weil complex can be used to compute $HC_{\delta}^*(\Ha)$, and that it arises naturally in the construction of characteristic maps associated to higher traces. We also include some examples. \\
\hspace*{.3in}{\bf\ \ \ \ \ \ \ \ \ \  Keywords}:\ \ Cyclic cohomology, Hopf algebras, $X$-complex, characteristic classes, Weil complex
\end{abstract}

\newpage

\tableofcontents

\newpage
  \hspace*{2in}{\bf INTRODUCTION}\\ 
\newline

\hspace*{.3in} In the transversal index theorem for foliation
(cohomological form), the characteristic classes involved are 
apriori cyclic cocycles on an algebra $\A$ associated to the
foliation. In their computation, A. Connes and
 H. Moscovici \cite{CoMo} have discovered that the
action of the operators appearing from the non-commutative index
formula can be organized in an action of a Hopf algebra $\Ha_{T}$ on
$\A$, and that the cyclic cocycles are made out just by combining the
 action with a certain invariant trace $\tau: \A\rmap {\bf C}$. In
 other words, they define a cyclic cohomology
 $HC_{\delta}^*(\Ha_{T})$, in such a way that the cyclic cocycles
 involved are in the target of a characteristic map $k_{\tau}:
 HC_{\delta}^*(\Ha_{T})\rmap HC^*(\A)$,  canonically associated to the
 pair $(\A,\tau)$. When computed, $HC_{\delta}^*(\Ha_{T})$ gives the
 Gelfand-Fuchs cohomology and the characteristic map $k_{\tau}$ is a
 non-commutative version of the classical \cite{Bo} characteristic map $k:
 H^*(WO_{q})\rmap H^*(M/\F)$ for codimension $q$ foliations $(M,\F)$
 (see \cite{BrNi, Cra} for the relation between $HC^*(\A)$ and $ H^*(M/\F)$). The definition of 
$HC_{\delta}^*(\Ha)$ given in \cite{CoMo}  applies to any Hopf algebra
 $\Ha$ endowed with a character $\delta$, satisfying certain
 conditions (see the end of \ref{problem}). In the context of $\Ha_{T}$ this provides a new beautiful
relation of cyclic cohomology with Gelfand-Fuchs cohomology,
 while, in general, it can be viewed
as a non-commutative extension of the Lie algebra cohomology (see
Theorem \ref{env}).\\
\hspace*{.3in}Our first goal is to show that Connes-Moscovici formulas can be used under the minimal requirement $S_{\delta}^2= Id$ (which answers a first question raised
in \cite{CoMo}), and to give a new
 definition/interpretation of $HC_{\delta}^*(\Ha)$ in the spirit of Cuntz-Quillen's formalism.\\
\hspace*{.3in} Independently, in the work of Gelfand and Smirnov 
on universal Chern-Simons classes, there are implicit relations with
cyclic cohomology (\cite{Gel}).
 Our second goal is to  make these connections explicit. This 
leads us to a noncommutative Weil complex $W(\Ha)$ associated to
a coalgebra, which extends the constructions from \cite{Gel, Chern}.\\
\hspace*{.3in}Our third goal is to show that $W(\Ha)$ is intimately related to
the cyclic cohomology $HC_{\delta}^*(\Ha)$ and to the construction (see \ref{daev}, \ref{daod}) 
of characteristic homomorphisms $k_{\tau}$ associated to invariant
higher traces $\tau$ (which is a second problem raised in
 \cite{CoMo}). The construction of $k_{\tau}$ is inspired by the
 construction of the usual Chern-Weil homomorphism (see e.g \cite{Duff}), and of the secondary characteristic classes for foliations
 (\cite{Bo}). \\
\hspace*{.3in}This work is strongly influenced by the Cuntz-Quillen
approach to cyclic cohomology (\cite{CQ3, CQ2, Qext} etc).\\
\hspace*{.3in} Here is an outline of the paper. In Section \ref{Prel} we bring together some 
basic results about characters $\delta$ and the associated twisted antipodes $S_{\delta}$ on
Hopf algebras. In section \ref{Itr} we present some basic terminology,
 describe the problem (see \ref{problem})  of defining
 $HC_{\delta}^*(\Ha)$, and explain why the case $S_{\delta}^2= Id$ is
 better behaved (see Proposition \ref{descend}). Under this
 requirement, we define a cyclic cohomology $HC_{\delta-inv}^*(R)$
 for any $\Ha$-algebra $R$, and we indicate how Cuntz-Quillen 
machinery can be adapted to this situation (Theorem \ref{homot}).
 The relevant information which is needed for the cyclic cohomology 
of Hopf algebras, only requires a small part of this machinery. This is
captured by a localized $X$-complex (denoted by $X_{\delta}(R)$); 
in Section \ref{SXoper}, after recalling the $X$-complex
 interpretation of $S$-operations (see \ref{Soper}), we introduce 
$X_{\delta}(R)$ (see \ref{Xinv}) and compute it in the case
where $R= T(V)$ is the tensor algebra of an $\Ha$-module $V$ (see Proposition \ref{bar}).\\
\hspace*{.3in} In Section \ref{Hocic} we prove that $HC_{\delta}^{*}(\Ha)$ can be defined under the minimal requirement $S_{\delta}^2=Id$ (see Proposition \ref{teorema}); also, starting with the question "which is the target of characteristic
 maps $k_{\tau}: HC_{\delta}^*(\Ha)\rmap HC^*(A)$, associated to pairs
 $(A,\tau)$ consisting of a $\Ha$-algebra $A$ and a $\delta$-invariant
 trace $\tau$ (see \ref{problem})?", we explain/interpret the definition of $HC_{\delta}^*(\Ha)$ in terms of localized $X$-complexes (see \ref{caract}, and Theorem \ref{teore}). This interpretation is the starting point in constructing the characteristic maps associated to higher traces (Section \ref{NCCW}). We also recall Connes-Moscovici's recent proposal to extend the definition to the non-unimodular case, and we show how our interpretations extend to this setting.\\
\hspace*{.3in} In Section \ref{Exs} we give some examples, including $\Ha= U_q(sl_2)$, and a
detailed computation of the fundamental example where $\Ha= U(\lie)$ is
the enveloping algebra of a Lie algebra $\lie$ (see Theorem \ref{env}).\\
\hspace*{.3in} In section \ref{NCW} we introduce the non-commutative
 Weil complex (by collecting together 'forms and curvatures' in a 
non-commutative way). We show that there are two relevant types of 
cocycles involved (which, in the case considered by Gelfand and 
Smirnov, correspond to Chern classes, and Chern-Simons classes,
respectively), we describe the Chern-Simons
transgression, and prove that it is an 
isomorphism between these two types of cohomologies (Theorem \ref{th1}). In connection
with cyclic cohomology, we show that the non-commutativity of the Weil
 complex naturally gives rise to an $S$-operator, and to cyclic
 bicomplexes computing our cohomologies (see \ref{sopt}-\ref{cordoi}).\\
\hspace*{.3in} In  section \ref{NCCW}
we come back to Hopf algebra actions, and higher traces, and we show how
the non-commutative Weil complex can be used to construct the characteristic map $k_{\tau}$ associated to higher traces (see \ref{daev},
\ref{daod}). To prove the compatibility with the $S$-operator (Theorems \ref{damm} and \ref{dam}), we also show that the truncations of the Weil complex can be expressed in terms of relative $X$-complexes (Theorem \ref{calcul}). When $\Ha= {\bf C}$, we re-obtain 
the cyclic cocycles (and their properties) described by Quillen
\cite{Qext}. In general, the truncated Weil complexes still
compute $HC^*_{\delta}(\Ha)$, as explained by Theorem \ref{comp}.
Section \ref{grea} is devoted to the proof of this theorem, and the
construction of characteristic maps associated to equivariant cycles.\\
\hspace*{.3in} {\bf Acknowledgments}: I would like to thank
I. Moerdijk for carefully reading the first (disastrous) version of
this work, and for making many valuable comments. I would also like to
thank J. Cuntz for his interest, and for inviting me to a workshop on
cyclic cohomology, where I attended a talk by A. Connes on this
subject. Also, I am thankful to V. Nistor for his help.

\section{Preliminaries on Hopf Algebras}
\label{Prel}

\hspace*{.3in} In this section we review some basic properties of Hopf algebras (see \cite{Sw})
and prove some useful formulas on twisted antipodes.\\
\hspace*{.3in} Let $\Ha$ be a Hopf algebra. As usual, denote by $S$ the antipode, by $\epsilon$
the counit, and by $\Delta(h)=\sum h_0\otimes h_1$ the coproduct. Recall some of the basic relations 
they satisfy:

\begin{equation} \sum \epsilon(h_0) h_1 =  \sum\epsilon(h_1) h_0 = h,\label{unu}\end{equation}
\begin{equation}\sum S(h_0) h_1 =  \sum h_0 S(h_1) = \epsilon(h)\cdot 1,\label{doi} \end{equation} 
\begin{equation}S(1)= 1, \epsilon(S(h))= \epsilon(h), \label{trei}\end{equation}
\begin{equation}S(gh)= S(h)S(g),\label{patru}\end{equation}  
\begin{equation}\Delta S(h)=\sum S(h_1)\otimes S(h_0).\label{cinci}\end{equation}
\hspace*{.3in}Throughout this paper, the notions of $\Ha$-module and $\Ha$-algebra have the 
usual meaning, with $\Ha$ viewed as an algebra. The tensor product $V\otimes W$ of two 
$\Ha$-modules is an $\Ha$-module with the diagonal action:
\begin{equation} h \ps v\otimes w\pd = \sum h_0\ps v\pd \otimes h_1 \ps w\pd. \label{diag}\end{equation}
\hspace*{.3in}A character on $\Ha$ is any non-zero algebra map $\delta :\Ha\rmap \bf{C}$. 
Characters will be used for 'localizing' modules: for
an $\Ha$-module $V$, define $V_{\delta}$ as the quotient of $V$ by the space of co-invariants 
(linear span of elements of type $h\ps v\pd -\delta \ps h\pd v$, with $h\in \Ha, v\in V$). In 
other words,
\[ V_{\delta} = \bf{C}_{\delta}\otimes _{\Ha} V,\]
where $\bf{C}_{\delta}=\bf{C}$ is viewed as an $\Ha$-module via $\delta$. Before looking at very
simple localizations (see \ref{local}), we need to discuss the 'twisted antipode' $S_{\delta}:\ = 
\delta \ast S$ associated to a character $\delta$ (recall that $\ast$ denotes the natural product
on the space of linear maps from the coalgebra $\Ha$ to the algebra $\Ha$, \cite{Sw}). 
Explicitly,
\[ S_{\delta}(h) = \sum \delta(h_0) S(h_1),\ \ \ \forall h\in \Ha .\]

\begin{lem} The following identities hold: 
\begin{equation} \sum S_{\delta}(h_0) h_1 = \delta(h)\cdot 1, \label{sapte}\end{equation}
\begin{equation} S_{\delta}(1)= 1, \epsilon(S_{\delta}(h))= \delta(h),\label{opt}
     \end{equation}
\begin{equation} \Delta S_{\delta}(h)= \sum S(h_1)\otimes S_{\delta}(h_0),\label{noua}
     \end{equation}
\begin{equation} S_{\delta}(gh)= S_{\delta}(h)S_{\delta}(g),\label{zece}\end{equation}
\begin{equation} \sum S^2(h_1) S_{\delta}(h_0)= \delta(h)\cdot 1 . \label{nou}\end{equation}
\end{lem}

\emph{proof:} These follow easily from the previous relations. For instance, the first relation follows from the definition of $S_{\delta}$, (\ref{doi}), and (\ref{unu}), respectively:
\[ \sum S_{\delta}(h_0)h_1 = \sum \delta(h_0) S(h_1)h_2 = \sum \delta(h_0)\epsilon(h_1)\cdot 1 = \delta( \sum h_0\epsilon(h_1))\cdot 1 = \delta(h)\cdot 1 .\]
The other relations are proved in a similar way.\ \ $\Boxe$\\


\begin{lem}. \label{imp} For any two $\Ha$-modules $V, W$:
\[  h\ps v\pd \otimes w \equiv v\otimes S_{\delta}(h)\ps w\pd \ \ \ \mbox{mod\ co-invariants}\]
\end{lem}
\emph{proof:} From the definition of $S_{\delta}$, $v\otimes S_{\delta}(h)w= \sum \delta(h_0) v\otimes S(h_1)w$, so, modulo coinvariants, it is $\sum h_0\ps v\pd\otimes h_1 S(h_2)\ps w\pd =$ $\sum \epsilon(h_1)h_0\ps v\pd \otimes w=$ $h\ps v
\pd \otimes w$, 
where for the last two equalities we have used (\ref{doi}) and (\ref{unu}), respectively.\ \ $\Boxe$ \\
\newline
\hspace*{.3in} It follows easily that:

\begin{cor}. \label{local} For any $\Ha$-module $V$, there is an isomorphism:
\[ (\Ha \otimes V)_{\delta} \equiv V,\ \ \ \  (h,v) \mapsto S_{\delta}(h)v.\]
\end{cor}


There is a well known way to recognize Hopf algebras with $S^2=Id$ (see \cite{Sw},
pp. 74). We extend this result to twisted antipodes:

\begin{lem}\label{twist}. For a character $\delta$, the following are equivalent:
\begin{enumerate}
\item[(i)] $ S^{2}_{\delta} = Id,$
\item[(ii)] $ \sum S_{\delta}(h_1)h_0 = \delta(h)\cdot 1, \ \forall h\in \Ha$
\end{enumerate} 
\end{lem}

\emph{proof:} The first implication follows by applying $S_{\delta}$ to (\ref{sapte}), using 
(\ref{zece}), and {\it (i)}. Now, assume {\it (ii)} holds. First, remark that $S\compose 
S_{\delta}= \delta $. Indeed,
\[ (S\ast (S\compose S_{\delta}))(h)= \sum S(h_0)S(S_{\delta}(h_1))=\sum S(S_{\delta}(h_1)h_0)= 
\delta (h)\cdot 1,\]
(where we have used the definition of $\ast$, (\ref{cinci}), and $(ii)$, respectively.) 
Multiplying this relation by $Id$ on the left, we get $S\compose S_{\delta}= Id\ast \delta$.
Using the definition of $S_{\delta}$, (\ref{noua}), and the previous relation, respectively,
\[ S_{\delta}^{2}(h)= \sum \delta(S_{\delta}(h)_0)S(S_{\delta}(h)_1)= \sum\delta(
S(h_1)S(S_{\delta}(h_0))= \sum \delta(S(h_2))h_0 \delta(h_1),\]
which is (use that $\delta$ is a character, and the basic relations again):
\[ \delta( \sum h_1 S(h_2))h_0 = \sum \delta(\epsilon(h_1))h_0= \sum \epsilon(h_1)h_0 = h.\ \ \Boxe \]


\section{Invariant traces}
\label{Itr}

In this section we present some basic terminology like invariant
traces, $\Ha$-algebras. For such an algebra $R$, the non-commutative
differential forms on $R$ can be localized, under the hypothesis
$S_{\delta}^2=Id$, and a cyclic cohomology $HC_{\delta-inv}^*(R)$
shows up. For completeness, we indicate how the Cuntz-Quillen
machinery \cite{CQ3} can be adapted to this context (see Theorem
\ref{homot}); this extends, in particular, the usual correspondence
(\cite{CQ3}) between $\delta$-invariant cyclic cocycles
and $\delta$-invariant higher traces (with equivariant linear splitting).

\begin{num}{\bf Flat algebras:}\label{flata}\emph{ Let $A$ be an
algebra, not necessarily unital. An action $\Ha\otimes A\rmap A$ of
$\Ha$ (viewed as an algebra) on $A$ is called flat (and say that $A$ is a $\Ha$-algebra) if:}\end{num}
\begin{equation}\label{flat} h(ab) = \sum h_0(a) h_1(b), \ \ \ \forall h\in \Ha,\ a, b\in A\end{equation}
The motivation for the terminology is that, in our interpretations (see \ref{caract}), it plays a role similar to the usual flat connections in geometry.

\begin{num}{\bf Invariant traces:}\label{TRAC}\emph{ Let $\Ha$ be a Hopf algebra endowed with a character $\delta$, and $A$ a $\Ha$-algebra. A trace $\tau: A \rmap \bf{C}$ is called $\delta$-invariant if:}
\[ \tau (ha) = \delta (h) \tau(a), \ \ \ \ \forall h\in \Ha ,\ a\in A .\]
\emph{If $\delta= \epsilon$ is the counit, we simply call $\tau$ invariant. \\
\hspace*{.3in}Recall \cite{Qext} that an even ($n$ dimensional) higher trace on an
algebra $R$ is given by an extension $0\rmap I\rmap
 L\rmap R \rmap 0$ and a trace on $L/I^{n+1}$,
 while an odd higher trace is given by an extension
as before, and an $I$- adic trace, i.e. a linear functional
 on $I^{n+1}$ vanishing on $[I^n, I]$. Starting with an extension 
of $\Ha$-algebras, and a $\delta$-invariant trace $\tau$, we talk 
about equivariant (or $\delta$-invariant) higher traces.}\end{num}

\begin{num}{\bf Examples:}\label{exemple}\emph{ If $\Ha= {\bf C}[\Gamma]$ 
is the group algebra of a discrete group 
$\Gamma$ (recall that $S(\gamma)= \gamma\otimes \gamma,
 \epsilon(\gamma)= 1$ if $\gamma = 1$, and $0$ otherwise),
 $\Ha$-algebras are precisely $\Gamma$-algebras.\\
\hspace*{.3in} If $G$ is a connected Lie group, $\lie$ its
 Lie algebra, and $\Ha= U(\lie)$ is the enveloping algebra,
 then $\Ha$-algebras are precisely infinitesimal $G$-algebras;
 that is, algebras $A$ endowed with linear maps (Lie derivatives)
$L_{v}: A\rmap A$, linear on $v\in \lie$, such that
 $L_{[v, w]}=L_vL_w-L_wL_v, L_v(ab)= L_v(a)b+ aL_v(b)$. 
If $\Delta: G\rmap {\bf C}$ is a character, it induces an
 infinitesimal character $\delta$ on $U(\lie)$:
 $\delta (v)= (\frac{d}{dt})_{t=0} \Delta(exp\ps tv\pd)$, $v\in \lie$.
 If $A$ is a topological (locally convex) $G$-algebra,
 $\delta$-invariance of traces on $A$ is equivalent to $\tau(ga)= 
\Delta(g)\tau(a), \ \forall g\in G, a\in A$.\\
\hspace*{.3in} Remark that, in general, the action of a
 Hopf algebra on itself is not flat. A basic example of 
flat action is the diagonal action of $\Ha$ on its tensor
 algebra $T\Ha= \oplus_{n\geq 1} \Ha^{\otimes n}$ (see (\ref{diag})).
 Another basic example is the algebra 
$\Omega^{*}(R)$ of noncommutative differential forms 
on a $\Ha$-algebra $R$. Recall that:}
\[ \Omega^n(R) = \tilde{R} \otimes R^{\otimes n}, \]
\emph{where $\tilde{R}$ is $R$ with a unit adjoined. Extending the action
 of $\Ha$ to $\tilde{R}$ by
$h\cdot 1:= \epsilon(h) 1$, we have an action of $\Ha$ on $\Omega^*(R)$ 
(the diagonal action). To check the flatness condition:
 $h(\omega\eta)= \sum  h_0(\omega)h_1(\eta),\ \forall \omega, \eta \, \in \Omega(R)$,
remark that one can formally reduce to the case where
 $\omega$ and $\eta$ are degree $1$ forms, in which case the computation is easy.\\
\hspace*{.3in}Recall also the usual operators $d, b, B, k$
 acting on $\Omega^*(R)$ (see  \cite{CQf}, paragraph $3$ of
 \cite{CQ3}):
 $d(a_0da_1 .\, .\, .\, . da_n)= da_0da_1 .\, .\, .\, . da_n,$ 
$b(\omega da)= (-1)^{deg\ps \omega\pd}[\omega, a]$, $ k= 1- (bd+ db)$,
 $B= (1+ k+ .\ .\ .\ + k^n)d$ on $\Omega^n(R)$.}\end{num}

\begin{num}{\bf The Problem:}\label{problem}\emph{ Let $\delta$ be a
 character on a Hopf algebra $\Ha$. The problem of defining a cyclic
 cohomology $'HC^*_{\delta}(\Ha)'$, should answer the question: which
 are the nontrivial cyclic cocycles on a $\Ha$-algebra $A$, arising
 from a $\delta$-invariant trace $\tau$, and the action of $\Ha$ on
 $A$. In particular, for any pair $(A, \tau)$ one should have an associated characteristic map:}
\[ k_{\tau}: HC^*_{\delta}(\Ha) \rmap HC^*(A)\ , \]
\emph{compatible with the $S$-operation on cyclic cohomology. There is
  a similar problem for invariant} higher \emph{traces.\\
\hspace*{.3in} In \cite{CoMo}, Connes and Moscovici have introduced $HC^*_{\delta}(\Ha)$ under the  hypothesis that there
is an algebra $\A$, endowed with 
an action of $\Ha$, and with a $\delta$-invariant \emph{faithful} trace 
$\tau: \A\rightarrow \bf{C}$. As pointed out in \cite{CoMo}, this requirement is quite strong; a more natural hypothesis would be the weaker condition $S_{\delta}^2= Id$.}\end{num}

\begin{prop}\label{descend}. 
If $S_{\delta}^2= Id$, then for any $\Ha$-algebra $R$, the operators $d, b, k, B$, acting on $\Omega^*(R)$ (see \ref{exemple}), descend to $\Omega^*(R)_{\delta}$.
\end{prop}

\emph{proof:} Since $d$ commutes with the action of $\Ha$, and $k, B$
(and all the other operators appearing in paragraph $3$ of \cite{CQ3})
are made out of $d, b$, it suffices to prove that, modulo co-invariants,
\[   b(h\cdot \eta) \equiv b(\delta(h) \eta), \ \ \ \forall h\in \Ha , \ \eta \in \Omega^*(R). \]
For $\eta = \omega da$, one has$(-1)^{|\omega|} b(h\cdot \omega da)= \sum h_0(\omega)h_1(a)-\sum h_1(a)h_0(\omega)$.\\
Using Lemma \ref{imp} and (\ref{sapte}), $\sum h_0(\omega)h_1(a)\equiv \sum \omega \cdot S_{\delta}(h_0)h_1 a = \delta(h)\omega a$.\\
Using Lemma \ref{imp}, and $(ii)$ of Lemma \ref{twist}, 
$\sum h_1(a)h_0(\omega)\equiv \sum a\cdot S_{\delta}(h_1)h_0\omega = \delta(h) a\omega$,
which ends the proof.\ \ $\Boxe$ \\

\begin{defin} ($S_{\delta}^2= Id$) Define the localized  cyclic
  cohomology $HC_{\delta-inv}^{*}(R)$ of $R$ as the cyclic 
cohomology of the mixed complex $\Omega^*(R)_{\delta}$. Similarly for
Hochschild and periodic cyclic cohomologies, and also for homology.
\end{defin}

\hspace*{.3in} This cohomology is not used in the next sections, but it fits very well in our discussion of higher traces. Recall that, via a certain notion of homotopy, higher traces correspond exactly to cyclic cocycles on $R$ (for
the precise relations, see pp. 417- 419 in \cite{CQ3}). Using
$HC^*_{\delta-inv}(R)$ instead of $HC^*(R)$, this relation extends to the
equivariant setting (provided one restricts to higher traces which
admit an equivariant linear splitting). The main ingredient is the following theorem which we include for completeness. It is analogous to one of the main results in \cite{CQ3} (Theorem 6.2). The notation $T(R)$ stands for the (non-unital) tensor algebra 
of $R$, and $I(R)$ is the kernel of the multiplication map $T(R)\rightarrow R$. Recall also
that if $M$ is a mixed complex \cite{Kassel}, $\theta M$ denotes the associated Hodge tower of $M$, which represents
the cyclic homology type of the mixed complex (for more details on the notations and terminology see \cite{CQ3}).

\begin{st} \label{homot}There is a homotopy equivalence of towers of super-complexes:
\[ \xx_{\delta}(TR, IR) \simeq \theta( \Omega^*(R)_{\delta}) .\]
\end{st}

\emph{proof:} The proof from \cite{CQ3} can be adapted. For this, one uses the fact that the projection $\Omega^*(R) \rmap \Omega^*(R)_{\delta}$ is compatible with all the structures
(with the operators, with the mixed complex structure). All the formulas we get 
for free, from \cite{CQ3}. The only thing we have to do is to take care of the action. 
For instance, in the computation of $\Omega^1(TR)_{\natural}$ (pp. $399- 401$ in \cite{CQ3}),
 the isomorphism $\Omega^1(TR)_{\natural}\cong \Omega^{-}(R)$  is not compatible with the 
action of $\Ha$, but, using the same technique as in  \ref{imp}, it descends to localizations  
(which means that we can use the natural (diagonal) action we have on $\Omega^{-}(R)$). 
With this in mind, the analogous of Lemma $5.4$ in \cite{CQ3} holds, that is, $\xx_{\delta}(TR, IR)$ 
can be identified (without regarding the differentials) with the tower 
$\theta(\Omega^*(R)_{\delta})$. Denote by $k_{\delta}$ the localization of $k$. 
The spectral decomposition with respect to $k_{\delta}$ is again a consequence of 
the corresponding property of $k$ (\cite{CQ3}, pp $389- 391$ and pp. $402-403$), and 
the two towers are homotopically concentrated on the nillspaces of $k_{\delta}$, 
corresponding to the eigenvalue $1$. Lemma $6.1$ of \cite{CQ3} identifies the two 
boundaries corresponding to this eigenvalues, which concludes the theorem.  \ \ $\Boxe$\\

\section{S-operations and X-complexes}
\label{SXoper}
In this section we recall Quillen's interpretation of a certain degree two cohomology operation ('$S$-operators') in terms of $X$-complexes, and describe a localized version (to be used in sections \ref{Hocic} and \ref{NCW}). As before, 
$\Ha$ is a Hopf algebra endowed with a character $\delta$ such that $S_{\delta}^2= Id$.
\begin{num}{\bf S-operations:}\label{Soper}\emph{ If $R$ is a DG
    algebra, denote by $R_{\natural}= R/[R, R]$ the complex obtained
    dividing out by the linear span of graded commutators. In examples
    like tensor algebras, the algebras considered by Gelfand, Smirnov
    etc (see \cite{Gel} and references therein), the noncommutative
    Weil complex of Section \ref{NCW}, and, in general  when $R$ is 'free', one encounters a very interesting degree two operation in the cohomology of $R_{\natural}$, $S: H^*(R_{\natural})\rmap H^{*+ 2}(R_{\natural})$. This phenomenon, due to the non-commutativity of $R$, has been very nicely explained by Quillen (\cite{Qext, Chern}). In general, for any algebra $R$, there is a sequence:}
\begin{equation}\label{short} 0\rmap R_{\natural} \stackrel{d}{\rmap} \Omega^1(R)_{\natural} \stackrel{b}{\rmap} R\stackrel{\natural}{\rmap} R_{\natural} \rmap 0, \end{equation}
\emph{Here $\Omega^1(R)_{\natural}= \Omega^1(R)/[\Omega^1(R), R]$, $b(xdy)= [x, y]$, and $\pi$ is the projection. In our graded setting, one uses graded commutators, and (\ref{short}) is a sequence of complexes. In general, it is exact in the right. When it is exact (and this happens in our examples), it can be viewed as an $Ext^2$ class, 
and induces a degree $2$ operator $S: H^{*}(R_{\natural})\rmap H^{*+2}(R_{\natural})$, explicitly described by the following diagram chasing (\cite{Chern}, pp. $120$). Given $\alpha \in H^k(R_{\natural})$, we represent it by a cocycle $c$, and use the exactness to solve successively the equations:}
\[ c= \natural(u) , \ \partial(u)= b(v) ,\ \partial(v)=  d(w) . \]
\emph{where $\partial$ stands for the vertical boundary. Then $S(\alpha)= [\, \natural(w)]\in H^{k+2}(R_{\natural})$.\\
\hspace*{.3in}Equivalently, pasting together (\ref{short}), we get a resolution, usually denoted by $X^{+}(R)$:}
\[ 0\rmap R_{\natural}\stackrel{d}{\rmap}  \Omega^1(R)_{\natural} \stackrel{b}{\rmap} R \stackrel{d}{\rmap} \Omega^1(R)_{\natural} \stackrel{b}{\rmap} R \rmap .\ .\ . \]
\emph{Emphasize that, when working with bicomplexes with anti-commuting differentials, one has to introduce a '$-$' sign for the even vertical boundaries (i.e. for those of $R$). So, one can use the cyclic bicomplex $X^+(R)$ to compute the cohomology of $R_{\natural}$, and then $S$ is simply the shift operator.\\
\hspace*{.3in} The $X$-complex of $R$ is simply the full version of $X^{+}(R)$, that is, the super-complex:}
\begin{eqnarray}\label{Xcom} X(R): \ \  \xymatrix{  R\ \ar@<-1ex>[r]_-{d} & \ \ \Omega^1(R)_{\natural} \ar@<-1ex>[l]_-{b}\ \ \ \ ,   }\end{eqnarray}
\emph{where $b(xdy)= [x, y], d(x)= dx.$ It is defined in general, for any algebra, and it can be viewed as the degree one level of the Hodge tower associated to $\Omega^*(R)$. In our graded setting, it is a cyclic bicomplex.}\end{num}

\begin{num}{\bf The localized $X$-complex:}\label{Xinv}\emph{ When $R= T\Ha$ is the tensor DG algebra of $\Ha$, then 
$T\Ha_{\natural}$ computes the cyclic cohomology of $\Ha$, viewed as a coalgebra (cf. Theorem \ref{quillen}), and our previous discussion describes the usual $S$-operation in cyclic cohomology. We need a similar construction for $T\Ha_{\natural, \delta}$.
Here, if $R$ is a DG algebra endowed with a flat action of $\Ha$ compatible with the differentials (a $\Ha$-DG algebra on short), $R_{\natural, \delta}:= R/[R,R]+({\rm coinvariants})$ denotes the complex obtained dividing out $R$ by the linear span of graded commutators and coinvariants (i.e. elements of type $h(x)- \delta(h)x$, with $h\in\Ha, x\in R$).\\
\hspace*{.3in} Since $S_{\delta}^2= Id$, we know (cf. Proposition \ref{descend}) that $b, d$ descend, and we define the localized $X$-complex of $R$ as the degree one level of the Hodge tower associated to $\Omega^*(R)_{\delta}$. In other words, this is simply the super-complex (a cyclic bicomplex in our graded setting):}
\[ X_{\delta}(R): \ \ \xymatrix{  R_{\delta}\ \ar@<-1ex>[r]_-{d} & \ \ \Omega^1(R)_{\natural, \delta} \ar@<-1ex>[l]_-{b} \ \ ,   } \]
\emph{where:}
\[ \Omega^1(R)_{\natural ,\delta}:\,=  \ \Omega^1(R)_{\delta}/b\Omega^2(R)_{\delta}= \ \Omega^1(R)/ [\Omega^1(R),R]+ ({\rm coinvariants})    ,\]
\emph{and the formulas for $b, d$ are similar to the ones for $X(R)$. There is one remark about the notation: $\Omega^1(R)_{\natural ,\delta}$ is not the localization
of $\Omega^1(R)_{\natural}$; in general, there is no natural action of $\Ha$ on it.}\end{num}

 
\begin{num}\label{tensor} {\bf Example.}
\emph{Before proceeding, let's look at a very important example: the 
(non-unital) tensor algebra $R= T(V)$ of an $\Ha$-module $V$. Adjoining a unit, one gets the
unital tensor algebra $\tilde{R}= \widetilde{T}(V)= \oplus_{n\geq 0} V^{\otimes n}$. The computation of
$X(R)$ was carried out in \cite{Qext}, Example 3.10. One knows that (\cite{CQ3}, pp. 395) $R= T(V), \ \ \Omega^1(R)_{\natural}=$ $V\oplus \widetilde{T} (V)= T(V)$,
and also the description of the boundaries: $d= \sum_{i=0}^{i=n} t^{i}, b= (t- 1)$ on
$V^{\otimes \ps n+1\pd}$, where $t$ is the backward-shift cyclic permutation. The second isomorphism is essentially due to the fact that, since $V$ generates $T(V)$, any element in $\Omega^1(T(V))$ can be written in the form $x d(v)y$, with $x, y\in \tilde{T(V)}, v\in V$ (see also the proof of the next proposition). To compute $X_{\delta}(R)$, one still has to compute its 
odd part. The final result is:}\end{num}

\begin{prop}\label{bar}. For $R= T(V)$:
\[ X_{\delta}^{0}(R)= T(V)_{\delta}, \ \ X_{\delta}^{1}(R)= T(V)_{\delta},\]
where the action of $\Ha$ on $T(V)$ is the usual (diagonal), and the boundaries have the 
same description as the boundaries of $X(R)$: they are $(t- 1)$, $N$ (which descend to the localization). The same holds when $V$ is a graded $\Ha$-module, provided we replace the the backward-shift cyclic permutation $t$ by its graded version. 
\end{prop}

\emph{proof:} One knows (\cite{CQ3},  pp. 395):
\[ \tilde{R}\otimes V\otimes \tilde{R} \tilde{\rmap} \Omega^1(R),\ \ x\otimes v\otimes y\mapsto
x(dv)y,\]
which, passing to commutators, gives (compare to \cite{CQ3}, pp. 395):
\[ R=V\otimes \tilde{R} \tilde{\rmap} \Omega^{1}(R)_{\natural}, \ \ v\otimes y\rmap \natural(dvy),\]
and the projection map $\natural: \Omega^1(R)\rmap \Omega^1(R)_{\natural}$ identifies with:
\[ \natural : \tilde{R}\otimes V\otimes \tilde{R} \rmap V\otimes \tilde{R}, \ \ 
x\otimes v\otimes y \rmap v\otimes yx .\]
So $X_{\delta}^1(R)$ is obtained from $T(V)$, dividing out by the linear subspace generated
by elements of type:
\[ \natural(h\cdot x\otimes v\otimes y- \delta(h)x\otimes v\otimes y)= 
\sum h^1(v)\otimes h^2(y)h^0(x) - \delta(h) v\otimes yx \in T(V).\]
Now, for $y= 1$, this means exactly that we have to divide out by coinvariants (of the diagonal 
action of $\Ha$ on $T(V)$). But this is all, because modulo these coinvariants we have (from 
Lemma \ref{imp}):
\[ \sum h_1(v)\otimes h_2(y)h_0(x)\equiv \sum v\otimes S_{\delta}(h_1)\cdot (h_2(y) h_0(x)), \]
while, from (\ref{noua}), (\ref{doi}) and $(ii)$ of Lemma \ref{imp}, (\ref{unu}):
\[ \sum S_{\delta}(h_1)\cdot (h_2(y) h_0(x))= \sum S(h_2)h_3(y) S_{\delta}(h_1)h_0(x)=
\sum \epsilon (h_1)\delta (h_0) yx= \delta (h) yx. \ \ \Boxe \]

\section{Cyclic Cohomology of Hopf Algebras}
\label{Hocic}

In this section we introduce the cyclic cohomology of Hopf algebras (endowed with a character $\delta$ as before). First we prove that Connes-Moscovici's formulas can be used under the minimal condition $S_{\delta}^2= Id$ (see \ref{teorema}). Next (see \ref{caract}) we present a second approach to defining $HC_{\delta}^*(\Ha)$ as the natural solution to our problem \ref{problem}. The two approaches coincide, which leads us to a $X$-complex interpretation of our cohomology (see Theorem \ref{teore}). This interpretation is also the starting point in dealing with higher traces (section \ref{NCCW}).\\
Let $\Ha$ be a Hopf algebra endowed with a character $\delta$. 
\begin{num}{\bf Cyclic cohomology of coalgebras:}\label{coalg}\emph{ Looking first just at the coalgebra structure 
of $\Ha$, one defines the cyclic cohomology of $\Ha$ by duality with the case of algebras. As in
\cite{CoMo}, we define the $\Lambda$-module (\cite{Co1}), denoted $\Ha^{\sharp}$, which is $\Ha^{\otimes\ps n+1\pd}$ in degree $n$,
whose co-degeneracies are:}\end{num}
\[ d^{\,i}(h^0,\, .\, .\, .\, , h^n) = \left \{ \begin{array}{ll}
         (h^0,\, .\, .\, .\, , h^{i-1}, \Delta h^{i}, h^{i+1},\, .\, .\, .\, , h^{n}) & \mbox{if $0 \leq i \leq n$} \\
         \sum (h^0_{\ps 1\pd}, h^1,\, .\, .\, .\,\, , h^{n}, h^{0}_{\ps 0\pd }) & \mbox{if $i=n+1$} 
                                      \end{array}
                             \right. .\]
and whose cyclic action is:
\[ t(h^0,\, .\, .\, .\, , h^n)= (h^1, h^2,\, .\, .\, .\, , h^n, h^0) .\]
Denote by $HC^*(\Ha)$ the corresponding cyclic cohomology, by $C^{*}_{\lambda}(\Ha)$ the cyclic 
complex, and by $CC^{*}(\Ha)$ the cyclic (upper plane) bicomplex (Quillen-Loday-Tsygan's) computing it.
Recall that  the DG tensor
algebra of $\Ha$, denoted $T(\Ha)$, is $\Ha^{\otimes n}$ in degrees $n\geq 1$ and $0$ otherwise, and has the differential $b\, '= \sum_{0}^{n} (-1)^{i} d^{\, i}$. The following proposition shows that the $S$-operator acting on $HC^*(\Ha)$  (apriori described by the shift on $CC^{*}(\Ha)$), is the $S$-operator described by an $X$-complex: 

\begin{prop}\label{quillen} Up to a shift on degrees, the cyclic bicomplex of $\Ha$, $CC^*(\Ha)$ coincides with the $X$-complex
of the DG algebra $T(\Ha)$, and the cyclic complex $C^{*}_{\lambda}(\Ha)$ is isomorphic to
$T(\Ha)_{\natural}$. This is true for any coalgebra.
\end{prop}

\emph{proof:} It follows from the computation in the proof of Lemma \ref{bar}, or by dualizing the analogous result for algebras (Theorem $4$ and Lemma $2.1$ of \cite{Qext}).

\ \ $\Boxe$ \\
\hspace*{.3in} Let us be more precise about the shifts. In a precise way, the proposition identifies $CC^{*}(\Ha)$ with the super-complex of complexes:
\[  .\ .\ .\ \rmap X^1(T\Ha)[-1]\rmap X^0(T\Ha)[-1]\rmap X^1(T\Ha)[-1]\rmap .\ .\ .\ \ \ ,\]
and gives an isomorphism:\ \ $C_{\lambda}^{*}(\Ha)\cong T(\Ha)_{\natural}[-1]$.

\begin{num}{\bf Cyclic cohomology of Hopf algebras:}\label{formule}\emph{ Localizing  the cyclic module $\Ha^{\sharp}$, we obtain a new object, denoted $\Ha_{\delta}^{\sharp}$. By Lemma \ref{local}, it is $\Ha^{\otimes n}$ in degree $n$, and the projection becomes:}
\[ \pi: \Ha^{\sharp} \rmap \Ha^{\sharp}_{\natural}, \ \pi(h^0, h^1, .\ .\ .\ , h^n)= S_{\delta}(h^0)\cdot (h^1, .\ .\ .\ , h^n),\]
\emph{where '$\cdot$' stands for the diagonal action of $\Ha$ (cf. Section $2$, (\ref{diag})).}\end{num} 
\hspace*{.3in}It is not true in general that the structure maps of $\Ha^{\sharp}$ descend to maps $d_{\delta}^i, s_{\delta}^i, t_{\delta}$ on $\Ha^{\sharp}_{\delta}$, but the compatibility with $\pi$ forces the following formulas, which make sense in general (compare to \cite{CoMo}, formulas $(37)-(40)$):
 \[ d_{\delta}^{\, i}(h^1,\, .\, .\, .\, , h^n) = \left \{ \begin{array}{lll}
         (1, h^1,\, .\, .\, .\, , h^n) &  \mbox{if $i=0$}\\
         (h^1, \, .\, .\, .\, , h^{i-1}, \Delta h^{i}, h^{i+1},\, .\, .\, .\, , h^{n}) & \mbox{if $1\leq i \leq n$} \\
         (h^1, .\, .\, .\, , h^n, 1) & \mbox{if $i=n+1$} 
                                      \end{array}
                             \right. .\]
\[ s^{\, i}_{\delta}(h^1,\, .\, .\, .\, , h^n) = (h^1, . . . , \epsilon(h^{i+1}),\, .\, .\, .\, , h^n) , \ \ 0\leq i\leq n-1 ,\]
\[ t_{\delta}(h^1,\, .\, .\, .\, , h^n) = S_{\delta}(h_1)\cdot (h^2,\, .\, .\, .\, , h^n, 1) \]
\hspace*{.3in} Apriori $\Ha_{\delta}^{\sharp}$ is just an $\infty$-cyclic \cite{FeTs}  module (in the sense that the cyclic relation $t_{\delta}^{n+1}= 1$ is not necessarily satisfied).  As pointed out by Connes and Moscovici, checking directly the cyclic relation $t_{\delta}^{n+1}= Id$ (which forces $S_{\delta}^2= Id$) is not completely trivial. They have proved it in \cite{CoMo} under the assumption mentioned in \ref{problem}.


\begin{prop}\label{teorema} Given a Hopf algebra $\Ha$ and a character $\delta$,  the previous formulas make $\Ha^{\sharp}_{\delta}$ into a cyclic module if and only if $S_{\delta}^2= Id$. More precisely:
\[ t_{\delta}^{n+1}(h^1, h^2, .\ .\ .\ , h^n)= (S_{\delta}^2(h^1), .\ .\ .\ , S_{\delta}^2(h^n))\ .\]
\end{prop}

\emph{Proof:} Dualizing the construction for algebras (see \cite{Cra, FeTs, Ni1}), to any coalgebra homomorphism $\theta: \Ha\rmap \Ha$ one associates a $\infty$-cyclic module $\Ha^{\sharp}(\theta)$. It is a slight modification of $\Ha^{\sharp}$ of \ref{coalg}, obtained by replacing $d^{n+1}, t$ in \ref{coalg} by:
\[ d^{\, n+1}(h^0,\, .\, .\, .\, , h^n) = \sum (h^0_{\ps 1\pd}, h^1,\, .\, .\, .\,\, , h^{n}, \theta(h^{0}_{\ps 0\pd })),\]
\[ t(h^0,\, .\, .\, .\, , h^n)= (h^1, h^2,\, .\, .\, .\, , h^n, \theta(h^0)) .\]
 Choosing $\theta:= S_{\delta}^2$, since $\pi: \Ha^{\sharp}(\theta) \rmap \Ha^{\sharp}_{\delta}$ is surjective, it suffices to show that $\pi$ 
 is compatible with the structure maps. The non-trivial formulas 
are $\pi d_{\theta}^{n+1}= d_{\theta}^{n+1}\pi$, 
$\pi t_{\theta}= t_{\theta}\pi$. We prove the last one.
 We need the following two relations which follow easily 
from (\ref{cinci}), (\ref{noua}):
\begin{equation} 
\Delta^{n-1} S_{\delta}(h)= \sum S(h_{\ps n\pd})\otimes .\ .\ .\
\otimes S(h_{\ps 2\pd})\otimes S_{\delta}(h_{\ps 1\pd}) , \label{star}
\end{equation}
\begin{equation}
 \Delta^{n-1} S_{\delta}S(h)= \sum S^2(h_{\ps 1\pd})\otimes .\ .\ .\
 \otimes S^2(h_{\ps n-1\pd})\otimes S_{\delta}S(h_{\ps n\pd}). \label{dstar}
\end{equation}
(where the sums are over $\Delta^{n-1}h= \sum h_{\ps 0\pd}\otimes . . . \otimes h_{\ps n\pd}$.) We have : 
\begin{eqnarray*}
\lefteqn{t_{\delta} \pi (h^0, h^1, .\ .\ .\ , h^n)=}\\
& & = \sum t_{\delta}(( S(h^0_{\ps n\pd}), .\ .\ .\  
, S(h^0_{\ps 2\pd}), S_{\delta}(h^0_{\ps 1\pd}))\star (h^1, .\ .\ .\ , h^n))\\
& & = \sum S_{\delta}(h^1) S_{\delta}S(h^0_{\ps n\pd}) \cdot
 (S(h^0_{\ps n-1\pd}), .\ .\ .\  , S(h^0_{\ps 2\pd}),
 S_{\delta}(h^0_{\ps 1\pd}), 1)\star (h^2, .\ .\ .\ , h^n, 1)\ .
\end{eqnarray*}
where $\star$ stands for the componentwise 'product' on $\Ha^{\otimes
  n}$. We want to prove it equals to
 $\pi t_{\theta}( h^0, h^1, . . . , h^n)=$ 
$S_{\delta}(h^1)\cdot (h^2, . . . , h^n, S^2_{\delta}(h^0))=$
 $S_{\delta}(h^1)\cdot (1, . . . , 1, S_{\delta}^2(h^0))\star$
 $(h^2, . . .  , h^n, 1)$, so it suffices to show that 
for any $h^0= h\in \Ha$:
\begin{equation}
\sum S_{\delta}S(h_{\ps n\pd})\cdot (S(h_{\ps n-1\pd}), .\ .\ .\ ,
S(h_{\ps 2\pd}), S_{\delta}(h_{\ps 1\pd}), 1)= 
(1, .\ .\ .\ , S_{\delta}^2(h)) \label{dedem}. 
\end{equation}
Using (\ref{dstar}), the left hand side is:
\begin{equation} \sum (S^2(h_{\ps n\pd})S(h_{\ps n-1\pd}), .\ .\ .\ , S^2(h_{\ps 2n-3\pd})S(h_{\ps 2\pd}), S^2(h_{\ps 2n-2\pd})S(h_{\ps 1\pd}), S_{\delta}S(h_{\ps 2n-1\pd})) \end{equation}
Using successively (\ref{nou}) for $\delta= \epsilon$, (\ref{unu}), and the coassociativity of $\Delta$, this is:
\begin{eqnarray*}
\lefteqn{\hspace*{.5in}\sum (1, . . . , 1, \epsilon(h_{\ps 2\pd}), S^2(h_{\ps 3\pd})S_{\delta}(h_{\ps 1\pd}), S_{\delta}S(h_{\ps 4\pd}))=}\\
&  = \sum (1, . . . , 1, S^2(h_{\ps 2\pd})S_{\delta}(h_{\ps 1\pd}), S_{\delta}S(h_{\ps 3\pd}))
 =\sum (1, . . . , 1, \delta(h_{\ps 1\pd}), S_{\delta}S(h_{\ps 2\pd}))\ , 
\end{eqnarray*}
by (\ref{nou}). Since $\sum \delta(h_{\ps 1\pd}) S_{\delta}S(h_{\ps 2\pd})=$ $S_{\delta}(\sum \delta(h_{\ps 1\pd})S(h_{\ps 2\pd}))=$ $S_{\delta}^2(h)$, we obtain the right hand side of (\ref{dedem}).
\ \ $\Boxe$

\begin{defin}\label{caractz} If $S_{\delta}^{2}= Id$, define  $HC^{*}_{\delta}(\Ha)$ as the
  cohomology defined by the cyclic module $\Ha^{\sharp}_{\delta}$; denote by $C^{*}_{\lambda ,\delta }(\Ha )$ the associated cyclic complex, and by $CC^{*}_{\delta}(\Ha)$ the associated cyclic bicomplex.
\end{defin}

In connection to our problem \ref{problem}, to any $\delta$-invariant trace $\tau$ on a $\Ha$-algebra $A$ one associates a characteristic map
$k^{\tau}: HC_{\delta}^{*}(\Ha) \rmap HC^{*}(A)$ ,
\begin{equation}
\label{reftau}
k^{\tau}(h_1, .\, .\, .\, , h_n)(a_0, a_1, .\, .\, .\, , a_n)= \tau(
a_0 h_1(a_1) .\, .\, .\,  h_n(a_n)),
\end{equation}
which is compatible with the $S$-operator (since it exists at the level of cyclic modules). Next we interpret/motivate this characteristic map, as well as the cohomology under discussion.

\begin{num}\label{caract}{\bf The (localized) characteristic map:}\emph{ Let $A$ be a $\Ha$-algebra, and let  $\tau :A\rmap 
\bf{C}$ be a trace on $A$. There is an obvious map induced in cyclic cohomology (which uses just the coalgebra structure of $\Ha$):}\end{num}
\begin{equation} \label{taau}
 \gamma^{\tau}: HC^{*}(\Ha)\rmap HC^{*}(A), \ (h^0,\, .\, .\, .\, , h^n) \rmap 
\gamma (h^0,\, .\, .\, .\, , h^n),
\end{equation}
\[ \gamma (h^0,\, .\, .\, .\, , h^n)(a_0,\, .\, .\, .\, , a_n)= \tau (h^0(a_0)\cdot\, .\, .\, .\, \cdot h^n(a_n)).\]
\hspace*{.3in}In order to find the relevant complexes in the case of invariant traces, we give a different interpretation of this simple map. We can view the action of $\Ha$ on $A$, as a linear
map:
\[ \gamma_{0}: \Ha \rmap Hom(B(A), A)^{1}= Hom_{lin}(A, A) \]
where $B(A)$ is the (DG) bar coalgebra of $A$. Recall that $B(A)$ is $A^{\otimes n}$ in degrees $n \geq 1$ and $0$ otherwise, with the coproduct:
\[ \Delta(a_1\otimes a_2 \otimes .\, .\, .\, \otimes a_n)= \sum_{i=1}^{n-1} (a_1\otimes .\, .\, .\, \otimes a_i)\otimes (a_{i+1}\otimes .\, .\, .\, \otimes a_n),\]
and with the usual $b^{\, '}$ boundary as differential. Then $Hom(B(A), A)$ is naturally a
DG algebra (see \cite{Qext}), with the product: $\phi *\psi:\ = m\compose (\phi\otimes\psi)\compose\Delta$ ($m$ stands for the multiplication on $A$). Explicitly, for $\phi, \psi \in Hom(B(A), A)$ of degrees $p$ and $q$, respectively,
\[ (\phi *\psi)(a_1, .\, .\, .\, , a_{p+q})= (-1)^{pq} \phi(a_1, .\, .\, .\, , a_{p}) \psi(a_{p+1}, .\, .\, .\, , a_{p+q})\ ,\]
\hspace*{.3in}The map $\gamma_{0}$ uniquely extends to a DG algebra map: 
\begin{equation}\label{buc1} \tilde{\gamma}: T(\Ha) \rmap Hom(B(A), A). \end{equation}
\hspace*{.3in}This can be viewed as a characteristic map for the flat action (see Proposition \ref{univ}).
Recall also (\cite{Qext}) that the norm map $N$ can be viewed as a closed cotrace $N: C_{*}^{\lambda}(A)[1]\rightarrow B(A)$ on the DG coalgebra $B(A)$, that is, $N$ is a chain map with the property that $\Delta\compose N= \sigma\compose\Delta\compose N$, where $\sigma$ is the graded twist $x\otimes y\mapsto
 (-1)^{deg\ps x\pd deg\ps y\pd} y\otimes x$. A formal property of this is that, composing with $N$ and $\tau$, we have an induced trace:
\begin{equation}\label{buc2} \tau_{\natural}: Hom(B(A), A)\rmap C_{\lambda}^{*}(A)[1],\ \ \tau_{\natural}(\phi)= \tau\compose\phi\compose N .\end{equation}
Composing with $\tilde{\gamma}$, we get a trace on the tensor algebra:
\begin{equation}\label{buc3} \tilde{\gamma}^{\tau}: T(\Ha)\rmap C_{\lambda}^{*}(A)[1], \end{equation}
and then a chain map:
\begin{equation}\label{buc4} \gamma^{\tau}: T(\Ha)_{\natural} \rmap C_{\lambda}^{*}(A)[1]. \end{equation}
\hspace*{.3in}Via Proposition \ref{quillen}, it induces (\ref{taau}) in cohomology.\\
\hspace*{.3in}Let's now start to use the Hopf algebra structure of $\Ha$, and the character $\delta$. First of all remark that the map $\tilde{\gamma}$ is $\Ha$-invariant, where the action of $\Ha$ on the right hand side of (\ref{buc1}) comes from the action on $A$: $(h\cdot\phi)(a)= h\phi(a),\ \forall\ a\in B(A)$. To check the invariance condition: $\tilde{\gamma}(hx)= h\tilde{\gamma}(x),\ \forall\ x\in T(\Ha)$, remark that the flatness of the action reduces the checking to the case where $x\in \Ha= T(\Ha)^{1}$, and that is obvious. Secondly, remark that if the trace $\tau$ is $\delta$-invariant, then so is (\ref{buc2}). In conclusion, $\tilde{\gamma}^{\tau}$ in (\ref{buc3}) is an invariant trace on the tensor algebra, so our map (\ref{buc4}) descends to a chain map:
\[ \gamma^{\tau}_{\delta}: T(\Ha)_{\natural, \delta} \rmap C_{\lambda}^{*}(A)[1],\]
So $H^*(TH_{\natural, \delta})$ naturally appears as the solution of our problem \ref{problem};  also, using the localized $X$-complex  $X_{\delta}(T\Ha)$ (see \ref{Xinv}), we have a short exact sequence:
\begin{equation}\label{sir}
 0\rmap T\Ha_{\natural, \delta} \stackrel{N}{\rmap} T\Ha_{(\delta)}
 \stackrel{1-t}{\rmap} T\Ha_{(\delta)} \rmap T\Ha_{\natural,
 \delta}\rmap 0
\end{equation}
describing an $S$-operation (cf \ref{Soper}) in our cohomology $H^*(TH_{\natural, \delta})$. These new objects are related to \ref{caractz} by the following (compare to Proposition \ref{quillen}):

\begin{st}\label{teore} Given a Hopf algebra $\Ha$ and a character $\delta$ such that $S_{\delta}^2= Id$, one has isomorphisms:
 \[C^{*}_{\lambda, \delta}(\Ha)\cong TH_{\natural, \delta}, \ \  \ CC^{*}_{\delta}(\Ha)\cong X_{\delta}(T\Ha), \]
up to the same degree shift as in Proposition \ref{quillen}. 
\end{st}

\emph{proof:} We have seen in Proposition \ref{bar}:
\[ \Omega^1(T\Ha)_{\natural} \cong T\Ha,\ \ \ \Omega^1(T\Ha)_{\natural, \delta} \cong (T\Ha)_{\delta} .\]
The first isomorphism is the one which gives the identification $X(T\Ha)\cong CC^*(\Ha)$ of Proposition \ref{quillen}. The second isomorphism, combined with the isomorphism (cf. Lemma \ref{local}):
\[ (T\Ha)_{\delta}^{n+1} \cong \Ha^{\otimes n}, \ [h_0\otimes h_1\otimes .\ .\ .\ \otimes h_n]\mapsto S_{\delta}(h_0)\cdot (h_1\otimes .\ .\ .\ \otimes h_n), \]
(with the inverse $ h_1\otimes .\ .\ .\ \otimes h_n \mapsto [1\otimes
h_1\otimes .\ .\ .\ \otimes h_n]$\,), gives the identification
$X_{\delta}(T\Ha)\cong CC_{\delta}^*(\Ha)$. \ \ $\Boxe$

\begin{num}\label{general}{\bf The uni-modular case:}\emph{ Motivated by examples like quantum groups, compact matrix groups \cite{Wo} and their duals, Connes and Moscovici have recently proposed \cite{unimod}  an extension of $HC_{\delta}^*(\Ha)$ to the more general case where $S_{\delta}$ is not necessarily involutive, but there exists an invertible group-like element $\sigma\in\Ha$ such that:
\begin{equation}\label{modular}  S_{\delta}^2(h)=\sigma h\sigma^{-1} \ \ \forall \ h\in \Ha, \ \ \delta(\sigma)=1 .
\end{equation}
In the terminology of \cite{unimod}, one says that $(\delta, \sigma)$ is a modular pair. For any such pair $(\delta, \sigma)$, one defines a cyclic module $\Ha^{\sharp}_{\delta,\sigma}$ by the same formulas as in \ref{formule} except for:
\[ d_{\delta,\sigma}^{\,n+1}(h^1,\, .\, .\, .\, , h^n)= (h^1,\, .\, .\, .\, , h^n, \sigma)\ ,\]
\[ t_{\delta,\sigma}(h^1,\, .\, .\, .\, , h^n)= S_{\delta}(h_1)(h^2,\, .\, .\, .\, , h^n, \sigma)\ .\]
Let $C_{\lambda, \delta, \sigma}^*(\Ha)$, $CC_{\delta, \sigma}^*(\Ha)$ be the associated cyclic complex, and cyclic bicomplex, respectively. The resulting cohomology is denoted by $HC_{\delta, \sigma}^*(\Ha)$, and appears as the target of characteristic maps associated to pairs $(A,\tau)$ with $\tau$ a $\delta$-invariant $\sigma$-trace (i.e. $\tau(ab)= \tau(b\sigma(a))$). 
}
\end{num}
Our interpretations extend to this setting. For any $\Ha$-algebra $R$, we define the following localized complex:
\[ X_{\delta, \sigma}(R): \ \ \xymatrix{  R_{\delta}\ \ar@<-1ex>[r]_-{d} & \ \ \Omega^1(R)_{\natural, \delta} \ar@<-1ex>[l]_-{b_{\sigma}} \ \ ,   } \]
where, this time, $b_{\sigma}(dx y)= -[x, y]_{\sigma}$, where $[x, y]_{\sigma}$ is the twisted commutator $xy- y\sigma(x)$, and $\Omega^{1}(R)_{\natural, \delta}$ is the quotient of $\Omega^{1}(R)$ by the subspace linearly spanned by coinvariants and twisted commutators $[x, \omega]_{\sigma}$ ($x\in R$, $\omega\in \Omega^{1}(R)$). Similarly one defines $R_{\natural, \delta}$, which fits into a sequence (exact on the right):
\[ 0\rmap R_{\natural, \delta} \stackrel{d}{\rmap} \Omega^1(R)_{\natural, \delta} \stackrel{b_{\sigma}}{\rmap} R_{\delta}\stackrel{\natural}{\rmap} R_{\natural, \delta} \rmap 0 .\]
The construction applies also to the graded case. Clearly, when $\sigma= 1$, we get the localized $X$-complex $X_{\delta}(R)$ defined in \ref{Xinv}.

\begin{st} Let $(\delta, \sigma)$ be as before. Then, for any $\Ha$ (DG) algebra $R$, $X_{\delta, \sigma}(R)$ is a well defined complex. For $R= T\Ha$:
\[ T\Ha_{\natural, \delta}\cong C_{\lambda, \delta, \sigma}^*(\Ha), \ \ \ , X_{\delta, \sigma}(T\Ha)\cong CC_{\delta, \sigma}^*(\Ha) \ ,\]
up to the same degree shift as in Proposition \ref{quillen}. 
\end{st}
\emph{proof:} The first part follows from the fact that $d: R\rmap \Omega^1(R)$, and $b_{\sigma}: \Omega^1(R)$ map coinvariants into coinvariants (with the same proof as for \ref{descend}), $b_{\sigma}$ kills the twisted commutators, $b_{\sigma}d= 0$ modulo coinvariants (straightforward), and $db_{\sigma}= 0$ in $\Omega^1(R)_{\natural, \delta}$. The last assertion follows from $\delta(\sigma)= 1$, and the relation:
\[ db_{\sigma}(dx y)= [x, dy]_{\sigma}- [\sigma^{-1}(y), dx]_{\sigma}+ (\sigma^{-1}(\omega)- \omega),\]
where $\omega= yd(\sigma(x))$. 
The second part is a straightforward extension of \ref{bar}, \ref{teore}. \ \ $\Boxe$ \\

One can also extend our interpretations \ref{caract} of the characteristic map.

\section{Some Examples}
\label{Exs}

In this section we compute the cohomology under discussion in several examples. Unless specified, $(\delta, \sigma)$ is a pair consisting of a character, and an invertible group-like element as in \ref{general} (i.e. satisfying $S_{\delta}^2(h)=\sigma h\sigma^{-1}$). In most of our examples, $\sigma=1$. \\
As a technical tool, let's remark that the complex computing $HH_{\delta,\sigma}^*(\Ha)$ depends just on the coalgebra structure of $\Ha$, and the group like elements $1,\sigma\in \Ha$. More precisely, denoting by ${\bf C}_{\sigma}$ the (left/right) one-dimensional $\Ha$ comodule induced by the group-like element $\sigma$, and by ${\bf C}$ the one corresponding to $\sigma= 1$, we have:


\begin{lem}\label{cotor} There are isomorphisms:
\[ HH_{\delta, \sigma}^{*}(\Ha) \cong Cotor_{\Ha}^*(\bf{C}, \bf{C}_{\sigma})\ ,\]
\end{lem}

\emph{proof:} For any group-like element $\sigma$ one has a standard resolution ${\bf C}_{\sigma}\stackrel{\sigma}{\rmap} B(\Ha, \bf{C}_{\sigma})$ of ${\bf C}_{\sigma}$ by (free) left $\Ha$ comodules. Here $B(\Ha, \bf{C}_{\eta})$ is $\Ha^{\otimes\ps n+1\pd}$ in degree $n$, and has the boundary:
\begin{equation}\label{dprim} d_{\sigma}{\, '}(h^0,\ .\ .\ .\ , h^n)= \sum_{i=0}^{n}(-1)^i (h^0,\ .\ .\ .\ , \Delta(h^i), \ .\ .\ .\ , h^n)+ (-1)^{n+1} (h^0,\ .\ .\ .\ , h^n, \sigma) .\end{equation}
Hence $Cotor_{\Ha}(\bf{C}, \bf{C}_{\sigma})$ is computed by the chain complex $\bf{C}\, \Boxe_{\Ha} B(\Ha, \bf{C}_{\sigma})$, that is, by the Hochschild complex of $\Ha^{\sharp}_{\delta}$.  \ \ \  $\Boxe$ \\

\begin{num}{\bf Example (group-algebras):}\emph{ If $\Ha= {\bf C}[\Gamma]$ is the group algebra of a discrete group $\Gamma$ (see \ref{exemple}), we have:
\[ HP_{\epsilon}^0({\bf C}[\Gamma]) \cong {\bf C}, \ \ HP_{\epsilon}^1({\bf C}[\Gamma]) \cong 0\]
($\epsilon=$ the counit, $\sigma= 1$).}
\end{num}

\emph{proof:}
We have the following periodic resolution $I^{*}$ of ${\bf C}$ by free ${\bf C}[\Gamma]$-comodules:
\[ 0\rmap {\bf C} \stackrel{\eta}{\rmap} {\bf C}[\Gamma]  \stackrel{\alpha}{\rmap} {\bf C}[\Gamma] \stackrel{\beta}{\rmap} {\bf C}[\Gamma]
\stackrel{\alpha}{\rmap} {\bf C}[\Gamma] \rmap .\ .\ .\ \]
where $\eta(1)=1$, $\alpha(g)= g$ for $g\neq 1$ and
 $\alpha(1)=0$, $\beta(g)= 0$ for $g\neq 1$  and $\beta(1)=1$. 
Hence $HH_{\epsilon}(\Ha)= Cotor_{\Ha}({\bf C}, {\bf C})$
 is computed by ${\bf C}\Boxe_{\Ha}I^{*}$, that is, by $0\rmap {\bf C}
\stackrel{0}{\rmap} {\bf C} \stackrel{id}{\rmap} {\bf C}\stackrel{0}{\rmap} {\bf C}
\stackrel{id}{\rmap} .\ .\ .\ $. So $HH_{\epsilon}^*(\Ha)= {\bf C}$
 if $n=0$ and $0$ otherwise, and the statement follows from the SBI sequence. \ \ \  $\Boxe$ \\

\begin{num}{\bf  Example (algebras with Haar integrals):}\emph{ Recall that a left Haar integral for the Hopf algebra $\Ha$ is a linear map $\tau:\Ha\rmap {\bf C}$ with the property $\tau(1)= 1$, $\sum \tau(h_0)h_1= \tau(h)\cdot 1$ for all $h\in \Ha$. Basic Hopf algebras which admit Haar integral are: finite dimensional Hopf algebras (by 5.1.6 of \cite{Sw}), group-algebras, algebras of smooth functions on a compact quantum group $G$ (by the fundamental Theorem 4.2 of \cite{Wo}). We recall that in the case of compact matrix groups there is a preferred choice of the character $\delta$, namely the modular character $f_{-1}$ of Theorem 5.6 \cite{Wo}. One has:
}
\end{num}

\begin{prop} If $\Ha$ admits a left Haar measure then:
\[ HP_{\delta}^0({\bf C}[\Gamma]) \cong {\bf C}, \ \ HP_{\delta}^1({\bf C}[\Gamma]) \cong 0\]
\end{prop}

\emph{proof:} Use the SBI sequence and the fact that the left integral $\tau$ induces a contraction $(h^0,\ .\ .\ .\ , h^n)$ $\mapsto$ $\tau(h^0) (h^1,\ .\ .\ .\ , h^n)$ of the Hochschild complex.  \ \ \  $\Boxe$ \\

\begin{num}{\bf Example (enveloping algebras):}\emph{ Following \cite{CoMo}
(Theorem $6. \ps i\pd$), We present now a detailed computation for the case where $\Ha=
U(\lie)$ is the  enveloping algebra of a Lie algebra $\lie$. Let $\delta$ be a character of $\lie$ (i. 
e. $\delta: \lie \rightarrow \bf{C}$ linear, with $\delta|_{[\lie,
\lie]}= 0$), and extend it to $U(\lie)$.  Denote by $\bf{C}_{\delta}$ the $\lie$ module $\bf{C}$
with the action induced by $\delta$. Since $S_{\delta}(v)=$ $-v +
\delta(v)$ for all $v \in \lie$, we are in the uni-modular case ($\sigma= 1$). 
The final result of our computation is:}\end{num}


\begin{st}\label{env} For any Lie  algebra $\lie$, and any $\delta \in \lie^*$:
\[HP_{\delta}^*(U(\lie)) \cong \bigoplus_{i\equiv * {\rm mod}\ 2} H_{i}(\lie ; \bf{C}_{\delta}) .\]
\end{st}

As a first step in the proof of \ref{env}, let's look at the symmetric
(Hopf) algebra $S(V)$ on a vector space $V$. Recall that the coproduct
is defined on generators by $\Delta(v)= v\otimes 1 + 1\otimes v,\
\forall v\in V$.


\begin{lem}\footnote{recently it has been pointed out to me that a
  version of this is due to P. Cartier, and appears also
  in \cite{Kass}, pp. 435- 442}\label{hoch} For any vector space $V$, the maps $A: \Lambda^n(V)\rmap S(V)^{\otimes n},\ v_1\wedge .\, .\, . \wedge v_n \mapsto (\sum_{\sigma} sign(\sigma) v_{\sigma\ps 1\pd} \otimes .\, .\, . \otimes v_{\sigma\ps n\pd})/n!$ induce isomorphisms:
\[ HH_{\delta}^{*}(S(V)) \cong \Lambda^*(V) .\]
\end{lem}

\emph{proof:} We will use a Koszul type resolution for the left $S(V)$ comodule $\bf{C}_{\delta}$. Let $e_1, . . . , e_k$ be a basis of $V$, and $\pi^i \in V^{*}$ the dual basis. The linear maps $\pi^i$ extend uniquely to derivations $\pi^i: S(V)\rmap S(V)$. Remark that each of the $\pi^i$'s are maps of left $S(V)$ comodules. Indeed, to check that $(1\otimes \pi^i)\compose \Delta= \Delta\compose \pi^i$, since both sides satisfy the Leibniz rule, it is enough to check it on the generators $e_i\in S(V)$, and that is easy. Consider now the co-augmented complex of left $S(V)$ comodules:
\[ 0\rmap \bf{C}_{\eta} \stackrel{\eta}{\rmap} S(V)\otimes \Lambda^0(V) \stackrel{d}{\rmap} S(V)\otimes \Lambda^1(V) \stackrel{d}{\rmap} .\ .\ .\ ,\]
with the boundary $d= \sum \pi^i \otimes e_i$, that is:
\[ d(x\otimes v_1\wedge .\, .\, . \wedge v_n)= \sum_{i=1}^{k} \pi^i(x)\otimes e_i\wedge v_1\wedge .\, .\, . \wedge v_n .\]
Point out that the definition does not depend on the choice of the basis, and it is dual to the Cartan boundary on the Weil complex of $V$, viewed as a  commutative Lie algebra. This also explains the exactness of the sequence. Alternatively, one can use a standard 'Koszul argument', or, even simpler, remark that $(S(V)\otimes\Lambda^*(V))\otimes (S(W)\otimes\Lambda^*(W))\equiv (S(V\oplus W)\otimes\Lambda^*(V\oplus W))$ as chain complexes (for any two vector spaces $V$ and $W$), which reduces the assertion to the case where $dim(V)= 1$. So we get a resolution $\bf{C}_{\eta} \rmap S(V)\otimes \Lambda^*(V)$ by free (hence injective) left $S(V)$ comodules. Then \ref{cotor} implies that $HH_{\delta}^*(S(V))$ is computed by $\bf{C}_{\eta}\Boxe_{S(V)} (S(V)\otimes \Lambda^*(V))$, that is, by $\Lambda^*(V)$ with the zero differential. This proves the second part of the theorem.\\

\hspace*{.3in}To show that the isomorphism is induced by $A$, we have to compare the previous resolution with the standard bar resolution $B(S(V), \bf{C}_{\eta})$ (see the proof of \ref{cotor}). We define a chain map of left $S(V)$ comodules:
\[ P: B(S(V), \bf{C}_{\eta}) \rmap S(V)\otimes \Lambda^*(V)\ ,\]
\[ P(x_0\otimes x_1\otimes .\, .\, . \otimes x_n)= x_0  \otimes pr(x_1)\wedge .\, .\, .  \wedge pr(x_n),\]
where $pr: S(V)\rmap V$ is the obvious projection map. We check now that it is a chain map, i.e.:
\[ dP(x_0\otimes x_1\otimes . . . \otimes x_n)= Pd(x_0\otimes x_1\otimes .\, .\, . \otimes x_n) .\]
First of all, we may assume $x_1, .\, .\, . , x_n \in V$ (otherwise, both terms are zero). The left hand side is then:
\[ \sum _{i=1}^{k}\pi^i(x_0)\otimes e_i\wedge x_1\wedge .\, .\, . \wedge x_n, \]
while the right hand side is:
\[ P(\Delta(x_0)\otimes x_1\otimes .\, .\, .  \otimes x_n)= (id\otimes pr)(\Delta\ps x_0\pd)\wedge x_1\wedge .\, .\, .  \wedge x_n\ .\]
So we are left with proving that:
\[ (id\otimes pr)\Delta(x) = \sum_{i=1}^{k} \pi^i(x)\otimes e_i , \ \ \forall \ x\in S(V) ,\]
and this can be checked directly on the linear basis $x= e_{i_1} .\, .\, . e_{i_n} \in S(V)$. In conclusion, $P$ is a chain map between our free resolutions  of $\bf{C}_{\eta}$ (in the category of left $S(V)$ comodules). By the usual homological algebra, the induced map $\bar{P}$ obtained after applying the functor $\bf{C}_{\eta}\Boxe_{S(V)} $-\, , induces isomorphism in cohomology. From the explicit formula:
\[ \bar{P}(x_1\otimes .\, .\, . \otimes x_n)= pr(x_1)\wedge .\, .\, . \wedge pr(x_n) ,\]
we see that $\bar{P}\compose A= Id$, so our isomorphism is induced by both $\bar{P}$ and $A$. \ \ \ $\Boxe$\\
\newline

\emph{proof of \ref{env}:} Consider the mixed complex \cite{Kassel}:
\[\Lambda : \xymatrix{ \Lambda^0(\lie) \ar@<-1ex>[r]_-{0} & \Lambda^1(\lie) \ar@<-1ex>[l]_-{d_{Lie}} \ar@<-1ex>[r]_-{0} & \Lambda^2(\lie) \ar@<-1ex>[l]_-{d_{Lie}} \ar@<-1ex>[r]_-{0} &\ .\ .\ .\ \ar@<-1ex>[l]_-{d_{Lie}} },\]
where $d_{Lie}$ stands for the usual boundary in the
Chevalley-Eilenberg complex computing $H_{*}(\lie)$. Denote by $\B$
the mixed complex associated to the cyclic module $\Ha^{\sharp}$, and
by $\B_{\delta}$ its localization, i.e. the mixed complex associated
to the cyclic module $\Ha^{\sharp}_{\delta}$ (so they are the mixed
complexes computing $HC^{*}(\Ha)$, and $HC_{\delta}^{*}(\Ha)$,
respectively). Here $\Ha= U(\lie)$. Let $\pi: \B\rightarrow
\B_{\delta}$ be the projection map, which, after our identifications
(see \ref{formule}), is degree-wise given by:
\[ \pi: \Ha^{\otimes\ps n+1\pd} \rmap \Ha^{\otimes n}, \ \pi(h_0\otimes .\, .\, . \otimes h_n)= S_{\delta}(h_0) \cdot (h_1\otimes .\, .\, . \otimes h_n).\]
Denote by $B$ and $B_{\delta}$ the usual (degree $(-1)$) '$B$- boundaries' of the two mixed complexes $\B, \B_{\delta}$.  Recall that $B= N\sigma_{-1} \tau$, where:
\[ \sigma_{-1}(h_0, \, .\, .\, .\, , h_n)= \epsilon(h_0) (h_1, \, .\, .\, .\, , h_n), \ \tau(h_0, \, .\, .\, . \, , h_n)= (-1)^n (h_1,  \, .\, .\, .\, ,  h_n, h_0)\ ,\]
and $N= 1+ \tau + .\, .\, . + \tau^n$ on $\Ha^{\otimes \ps n+1\pd}$.\\
\hspace*{.3in} We will show that $\Lambda$ and $\B_{\delta}$ are quasi-isomorphic mixed complexes (which easily implies the theorem), but for the computation we have to use the mixed complex $\B$, where explicit formulas are easier to write. We define the map:
\[ A: \Lambda^n(\lie) \rmap \Ha^{\otimes n},\ A(v_1\wedge .\, .\, . \wedge v_n)= (\sum_{\sigma} sign(\sigma) v_{\sigma\ps 1\pd} \otimes .\, .\, . \otimes v_{\sigma\ps n\pd})/ n! .\]
The fact that the (localized) Hochschild boundary depends just on the coalgebra structure of $U(\lie)$ and on the unit, which are preserved by the Poincare-Birkhoff-Witt Theorem (see e.g. \cite{We1}), together with Lemma \ref{hoch}, shows that $A$ is a quasi-isomorphism of mixed complexes, provided we prove its compatibility with the degree $(-1)$ boundaries, that is:
\begin{equation}\label{compat} B_{\delta}(A(x))= A(d_{Lie}(x)), \ \ \forall \ x= v_1\wedge \, .\, .\, . \wedge v_n \in \Lambda^n(\lie) .\end{equation} 
Using that $A(x)= \pi(y)$, where $y= (\sum sign(\sigma) 1\otimes v_{\sigma\ps 1\pd} \otimes .\, .\, . \otimes v_{\sigma\ps n\pd})/n!$,
 we have 
\begin{eqnarray*}
\lefteqn{B_{\delta}A(x)= \pi(B(y))= }\\
& & = \pi N\sigma_{-1} (\sum sign(\sigma)\,  (1\otimes v_{\sigma\ps 1\pd} \otimes .\, .\, . \otimes v_{\sigma\ps n\pd}- 
(-1)^nv_{\sigma\ps 1\pd} \otimes .\, .\, . \otimes v_{\sigma\ps n\pd}\otimes 1))/n!\\
& & = \pi(N(\sum sign(\sigma) v_{\sigma\ps 1\pd} \otimes .\, .\, . \otimes v_{\sigma\ps n\pd}))/n!\\
& & = \pi (\sum sign(\sigma) v_{\sigma\ps 1\pd} \otimes .\, .\, . \otimes v_{\sigma\ps n\pd})/(n-1)!\ .
\end{eqnarray*} 
But 
\[ \pi (v\otimes v_1\otimes .\, .\, . \otimes v_n) = \delta(v)  v_1\otimes .\, .\, . \otimes v_n - \sum_{i=1}^{n} v_1\otimes .\, .\, . \otimes v v_i\otimes .\, .\, . v_n\ ,\]
and, with these, it is straightforward to see that $B_{\delta}A(x)$ equals to:
\[A(\sum_{i=1}^n (-1)^{i+1} \delta(v_1) v_1\wedge .\, .\, . \wedge \widehat{v_i}\wedge .\, .\, . \wedge v_n + \sum_{i< j} (-1)^{i+ j} [v_i, v_j] \wedge v_1\wedge .\, .\, . \wedge \widehat{v_i}\wedge .\, .\, . \wedge \widehat{v_j}\wedge .\, .\, . \wedge v_n ),\]
i.e. with $A(d_{Lie}(x))$. \ \ \ $\Boxe$

\begin{num}{\bf Example (the quantum enveloping algebra of $sl_2$):}\emph{ We look now at the simplest example of a quantized envelopping algebra, namely $U_q(sl_2)$. As an algebra, it is generated by the symbols $E, F, K, K^{-1}$, subject to the relations: $KE= q^2EK, KF= q^{-2}FK, KK^{-1}= K^{-1}K= 1, [E, F]= (K- K^{-1})/(q- q^{-1})$. The co-algebra structure is given by:
\[ \Delta(K)= K\otimes K, \ \Delta(E)= 1\otimes E+ E\otimes K, \ \Delta(F)= K^{-1}\otimes F+ F\otimes 1,\]
\[ \epsilon(K)= 1,\  \epsilon(E)=  \epsilon(F)= 0 ,\]
while for the antipode: $S(K)= K^{-1},\, S(E)= -EK^{-1},\, S(F)= -KF$. One has $S^2(h)= KhK^{-1}$ for all $h\in U_q(sl_2)$, hence this is a first example with $\sigma\neq 1$.
}
\end{num}

\begin{prop} $HP^{0}(U_q(sl_2))= 0$, and $HP^{1}(U_q(sl_2))\cong {\bf C}^2$ with the generators represented by $E$ and $KF$.\\
(we have ommited the indices $\delta=\epsilon$, $\sigma=K$ from the notation)
\end{prop}

\emph{proof:} Denote $E= x$, $K= \sigma$, $KF= y$. Clearly $x$, and $y$ define cyclic cocycles. 
Using the SBI sequence, it suffices to prove the similar statement for Hochschild cohomology. We first prove that, for any $n$, 
\begin{equation}\label{negative} Cotor_{\Ha}({\bf C}, {\bf C}_{\sigma^a})= 0, \ \ \forall\ a\, < 0\ . \end{equation}
We use induction over $n$. It is obvious for $n= 0$; let's assume it is true for any $k< n$. Remark that, by the proof of \ref{cotor}, $Cotor({\bf C}_{\alpha}, {\bf C}_{\beta})$ ($\alpha, \beta$- group-like elements) is computed by the complex $B(\Ha; {\bf C}_{\alpha}, {\bf C}_{\beta})$, which is $\Ha^{\otimes n}$ in degree $n$ and has the boundary $u\mapsto (\alpha, u)- d_{\beta}{'}(u)$ (see (\ref{dprim})). Denote $B_a= B(\Ha; {\bf C}, {\bf C}_{\sigma^a})$, and $d_a$ its boundary.\\One has the following basis of $\Ha$:
\begin{equation}\label{basis} \{ x^my^k\sigma^p : m, k, p \ -{\rm integers}\ \ , m, k \geq 0\}\ . \end{equation}
Let '$\leq$' be the order $(m, k)\leq (m\,', k\,')$ iff $m\leq m\,'$, or $m= m\,'$ and $k\leq k\,'$. For any pair $(m, k)$ of positive integers, denote by $L_{m, k}$, and $L^{<}_{m, k}$ the subcomplexes of $B_a$ 
linearly spanned by elements of type $(x^iy^j, .\, .\, .\, )$,  with $(i, j)\leq (m, k)$, and $(i, j)< (m, k)$, respectively. \\
For the proof of (\ref{negative}), let $[z]\in Cotor^{n}({\bf C}, {\bf C}_{\sigma^a})$, represented by a cocycle $z\in B_a$. We claim that:
\begin{equation}\label{somem} \exists m, k\geq 0,\ \exists u\in L_{m, k}\ :\ \ [z]= [u] \ .\end{equation}
Indeed, defining $\tau:\Ha\rmap {\bf C}$ on the basis (\ref{basis}) by $\tau(1)= 1$ and $0$ otherwise, and $\theta= \tau\otimes Id_{\Ha}\otimes .\, .\, .\, \otimes Id_{\Ha}: B_a\rmap B_a$, we have $d_{a}\theta+ \theta d_{a}= Id-\phi$ where $\phi$ is identity on elements of type $(x^my^k, .\, .\, .\, )$, $(m, k)\neq (0, 0)$, and $0$ otherwise. Choose $(m, k)$ minimal such that $(\ref{somem})$ holds. Hence we find $v\in B_{a}^{n-1}$ such that:
\[ u\equiv (x^my^k, v)\  \ {\rm mod}\ \ \ L^{<}_{m, k}\ .\]
Assume first $(m, k)\neq (0, 0)$. Since $d_a(u)\equiv -(x^my^k,\sigma^{m+ k}, v)+ (x^my^k, d_{a}\,'(v))$ mod $L^{<}_{m, k}$, we must have $(\sigma^{m+k}, v)= d_{a}\,'(v)$, i.e. $v$ represents a $(n-1)$-cocycle in the standard complex computing $Cotor({\bf C}_{\sigma^{m+ k}}, {\bf C}_{\sigma^a})$. This complex is isomorphic (by the multiplication by $\sigma^{-m- k}$) to the standard complex computing $Cotor({\bf C}, {\bf C}_{\sigma^{a-m- k}})$, hence, by the induction hypothesis, $v= (\sigma^{m+ k}, w)- d_{a}\,'(w)$ for some $w$. 
Choosing $u\,'= u+ d_a(x^my^k, w)$, we then have $[z]= [u\,']$, and $u\,'\in L^{<}_{m, k}$, which contradicts the minimality of $(m, k)$.
We are left with the case $(m, k)= (0, 0)$, when, since $\phi(u)= 0$, one gets $[u]= 0$.\\
A completely similar argument shows that:
\[ Cotor^n({\bf C}, {\bf C})= 0\ \ \forall\ n\geq 1\ , \ \ Cotor^n({\bf C}, {\bf C}_{\sigma})= 0\ \ \forall\ n\geq 2\ ,\]
while clearly $Cotor^0({\bf C}, {\bf C}_{\sigma})= 0$. Let $[z]\in Cotor^1({\bf C}, {\bf C}_{\sigma})$. As above, we find a minimal $(m, k)$ such that (\ref{somem}) holds, and let $\eta \in {\bf C}$ such that $u\equiv \eta x^my^k$ mod $L^{<}_{m, k}$. Again, since $d_1(u)\equiv \eta( (x^my^k,\sigma^{m+ k})- (x^my^k,\sigma))$ mod $L^{<}_{m, k}$, we must have $\eta\sigma^{m+ k}= \eta\sigma$, hence $\eta= 0$, or $m+ k=1$. In other words, $Cotor^1({\bf C}, {\bf C}_{\sigma})$ $\cong {\bf C}^2$, with the generators $[x], [y]$. \ \ \ $\Boxe$

\section{A non-commutative  Weil complex}
\label{NCW}
In this section we introduce/describe a non-commutative Weil complex associated to a coalgebra, which extends/explains some results in \cite{Gel, Chern}, and will naturally appear in the construction of characteristic maps associated to higher traces (section \ref{NCCW}). We describe the relevant  cohomologies (analogues of Chern, Chern/Simons classes), and (using section \ref{SXoper}) the $S$-operators acting on them.\\

\hspace*{.3in}Let $\Ha$ be a coalgebra. Define its Weil algebra $W(\Ha)$ as the 
(non-commutative, non-unital) DG algebra freely generated by the symbols $h$ of degree $1$, $\omega({h})$ of degree
$2$, linear on $h\in \Ha$. The differential of $W(\Ha)$ is similar to
the $b^{\, '}$ differential of $T(\Ha)$ (see \ref{coalg}). It is denoted by $\del$, and is the unique derivation which acts on generators by:
\[ \del(h)= \omega_{h} - \sum h_0h_1 ,\]
\[ \del(\omega_{h})= \sum \omega_{h_0}h_{1}- \sum h_0\omega_{h_1} .\]

\begin{exam}\label{simple}\emph{ This algebra is intended to be a non-commutative analogue of the usual
Weil complex of a Lie algebra (see \cite{Duff}). Particular cases have been used in the study of universal Chern-Simons 
forms. When $\Ha = {\bf C}\rho$ (i.e. ${\bf C}$, with $1$ denoted by $\rho$), with $\Delta(\rho)=\rho\otimes \rho$, it is the algebra
introduced in \cite{Chern}; for $\Ha = {\bf C}\rho_1\bigoplus . . . \bigoplus {\bf C}\rho_n$ with
$\Delta(\rho_i)= \rho_i\otimes \rho_i$, we obtain one of the algebras studied on \cite{Gel}.}
\end{exam}

We discuss now its `universal property`. Given a DG algebra $\Omega^*$, and a linear map:
\[ \phi : \Ha \rmap \Omega^1 ,\]
define its curvature: \[ \omega_{\phi}: \Ha \rmap \Omega^2,\  \omega_{\phi}(h)= d\phi(h) + \sum \phi(h_0)\phi(h_1).\]
\hspace*{.3in}Alternatively, using the natural DG algebra structure of $Hom(\Ha ,\Omega^*)$,
\[ \omega_{\phi}:= d(\phi) - 1/2 [\phi ,\phi] \in Hom(\Ha ,\Omega)^2 .\]
\hspace*{.3in}There is a unique algebra homomorphism (the characteristic map of $\phi$):
\[ k(\phi): W(\Ha) \rmap \Omega^* ,\]
sending $h$ to $\phi(h)$ and $\omega_{h}$ to $\omega_{\phi}(h)$. \\
\hspace*{.3in}One can easily see that (compare with the usual Weil complex of a Lie algebra):

\begin{prop}\label{univ}. The previous construction induces a $1-1$ correspondence between linear maps 
$\phi : \Ha \rmap \Omega^1$ and DG algebra maps $k: W(\Ha) \rmap \Omega^*$. In particular, there is a $1-1$ correspondence between flat linear maps 
$\phi : \Ha \rmap \Omega^1$ (i. e. with the property that $\omega_{\phi}= 0$), and DG algebra maps $k:T(\Ha)\rmap \Omega^*$.
\end{prop}
An immediate consequence is that $W(\Ha)$ does not depend on the co-algebra structure of $\Ha$.
Actually one can see directly that $(W(\Ha), \del) \cong (W(\Ha), d)$, where $d$ is the 
derivation on $W(\Ha)$ defined on  generators by $d(h)= \omega_{h}, \ \ d(\omega_{h})= 0$ 
(i.e. the differential corresponding to $\Ha$ with the trivial co-product). An explicit isomorphism sends $h$ to
$h$ and $\omega_{h}$ to $\omega_{h} +\sum h_0h_1$. 
\begin{cor}\label{contract}. The Weil algebra $W(\Ha)$, and the complex $W(\Ha)_{\natural}$ are
 acyclic. \end{cor}

\begin{num} \label{bigrad}{\bf Extra-structure on $W(\Ha)$:}\emph{ Now we look at the extra-structure of $W(\Ha)$. First of all, denote by $I(\Ha)$ the ideal 
generated by the curvatures $\omega_{h}$. The powers of $I(\Ha)$, and the induced truncations
are denoted by:}
\[ I_{n}(\Ha):= I(\Ha)^{n+1},\ \ W_n(\Ha):= W(\Ha)/I(\Ha)^{n+1}.\]
\emph{Remark that $W_0(\Ha)= T(\Ha)$ is the tensor (DG) algebra of
  $\Ha$ (up to a minus sign in the boundary, which is irrelevant, and
  will be ignored). Dual to even higher traces, we introduce the complex:}
\[ W_n(\Ha)_{\natural}:=W_{n}(\Ha)/[W_{n}(\Ha), W_{n}(\Ha)]\] 
\emph{obtained dividing out the (graded commutators). In the terminology of \cite{Gel} (pp. 103), it is the space of 'cyclic words'. Dual to odd higher traces:}
\[ I_n(\Ha)_{\natural}:= I_n(\Ha)/[I(\Ha), I_{n-1}(\Ha)].\]\end{num}
\hspace*{.3in}It is interesting that all these complexes compute the same
cohomology (independent of $n$\, !), namely the cyclic cohomology of $\Ha$
viewed as a coalgebra. This is the content of Theorem \ref{th1},
Proposition \ref{quillen}, and Section \ref{grea}.\\
\hspace*{.3in}Secondly, we point out a bi-grading on $W(\Ha)$:
defining $\del_0$ such that $\del= \del_0 + d$, then $W(\Ha)$ has a
structure of bigraded differential algebra, with $deg(h)=(1, 0), deg(\omega_{h})= (1, 1)$.
Actually $W(\Ha)$ can be viewed as the tensor algebra of $\Ha^{\ps 1,0\pd}\bigoplus \Ha^{\ps 1,1\pd}$
 (two copies of $\Ha$ on the indicated bi-degrees). With this bi-grading, $q$ in
$W^{p, q}$ counts the number of curvatures. The  boundary $d$ increases $q$, while
$\del_0$ increases $p$.

\begin{num}{\bf Example:}\emph{ Let's have a closer look at $\Ha = {\bf C}\rho$ with $\Delta(\rho)= h\rho\otimes \rho$, for which the computations were carried out by D. Quillen \cite{Chern}, recalling the main features of our complexes:\\
\hspace*{.25in}$(1)$ $\omega^n$ are cocycles of $W(\Ha)_{\natural}$ (where $\omega= \omega_{h\rho}$). They are trivial in cohomology (cf. Corollary \ref{contract}).\\
\hspace*{.25in}$(2)$ the place where $\omega^n$ give non-trivial cohomology classes is $I_{m, \natural}$, with $m$ sufficiently large. \\
\hspace*{.25in}$(3)$ the cocycles $\omega^n$ (trivial in $W(\Ha)_{\natural})$ transgress to certain (Chern-Simons) classes. The natural complex in which these classes are non-trivial (in cohomology) is $W_m(\Ha)_{\natural}$.\\
\hspace*{.25in}$(4)$ there are striking 'suspensions' (by degree $2$ up) in the cohomology of all the complexes $W_n(\Ha)_{\natural}$, $I_n(\Ha)_{\natural}$, $\tilde{I}_{n} (\Ha)_{\natural}$.\\
\hspace*{.3in}Our intention is also to explain these phenomena (in our general setting).
}\end{num}

\begin{num}{\bf 'Chern-Simons contractions'}\label{path}. \emph{ Starting with two linear maps:}\end{num}
\[ \rho_0 , \rho_1 : \Ha \rmap \Omega^1 ,\] 
we form:
\[ t\rho_0 + (1-t)\rho_1: = \rho_0\otimes t + \rho_1\otimes (1-t) : \Ha \rmap (\Omega^*\otimes \Omega(1))^1,\]
where $\Omega(1)$ is the algebraic DeRham complex of the line: $\bf{C} [t]$ in degree $0$, and
$\bf{C} [t]dt$ in degree $1$, with the usual differential. Composing its characteristic map 
$W(\Ha) \rmap \Omega\otimes \Omega(1)$, with the degree $-1$ map 
$\Omega\otimes\Omega(1)\rmap\Omega$ coming from the integration map
$\int_{0}^{1} :\Omega(1)\rmap \bf{C}$ (emphasize that we use the graded tensor product, and the integration map has degree $-1$), we get a degree $-1$ chain map:
\[ k(\rho_0, \rho_1): W(\Ha) \rmap \Omega.\]
As usual, $[k(\rho_0, \rho_1), \del] = k(\rho_1) - k(\rho_0)$.\\ 
The particular case where $\Omega= W(\Ha), \rho_0= 0, \rho_1= Id_{\Ha}$ gives a contraction of $W(\Ha)$:
\[ H:= k(Id_{\Ha}, 0): W(\Ha)\rmap W(\Ha). \]
We point out that $H$ preserves commutators, and induces a chain map
\begin{equation}\label{chsi}
CS : I_{n}(\Ha)_{\natural}\rmap W_n(\Ha)_{\natural}[1] ,\ \ [x]\mapsto
[H(x)] ,
\end{equation}
to which we will refer as the Chern Simons map. The formulas for the
 contraction $H$ resemble the usual ones (\cite{Gel, Qext,
 Chern}). For instance, at the level of $W(\Ha)_{\natural}$, one has:
\begin{equation}\label{sim}
 H(\frac{\omega_{h}^{n+ 1}}{(n+1)!})=  \int_{0}^{1} 
\frac{1}{n!}h (t\omega_h + (t^2- t) \sum_{(h)} h_0h_1)^n dt 
\end{equation}


\begin{st}\label{th1} The Chern-Simons map is an isomorphism $H^*(I_{n}(\Ha)_{\natural}) \tilde{\rmap} H^{*- 1}(W_n(\Ha)_{\natural})$ (compatible with the $S$-operator described below). \end{st}

\emph{proof:} We consider the following slight modification of $I_n(\Ha)$:
\[ \tilde{I}_{n}(\Ha)_{\natural}:= I_n(\Ha)/I_{n}(\Ha)\cap [W(\Ha), W(\Ha)].\]
One has a short exact sequence:
\[ 0\rmap \tilde{I}_{n}(\Ha)_{\natural}\rmap W(\Ha)_{\natural}\rmap
W_n(\Ha)_{\natural}\rmap 0 , \]
and, using Corollary \ref{contract}, the boundary of the long 
exact sequence induced in cohomology gives an isomorphism 
$\tilde{\partial}: H^{*- 1}(W_n(\Ha)_{\natural}) \tilde{\rmap}
H^*(\tilde{I}_{n}(\Ha)_{\natural})$.
The same formula (\ref{chsi}) defines a chain map $\tilde{CS}:
\tilde{I}_{n}(\Ha)_{\natural}\rmap W_n(\Ha)_{\natural}[1]$, and one
can easily check that $\tilde{CS}\compose \tilde{\partial}= Id$. Now, since
$CS$ is the composition of $\tilde{CS}$ with the canonical projection
$I_n(\Ha)_{\natural}\twoheadrightarrow \tilde{I}_{n}(\Ha)_{\natural}$,
it suffices to show that the last map induces isomorphism in
cohomology. We prove this after describing the $S$- operator.  \ \ $\Boxe$


\begin{num}{\bf The $S$-operator:}\label{sopt}\emph{ The discussion in \ref{Soper} applies to the Weil complex $W(\Ha)$, explaining the 'suspensions' (by degree $2$ up) in the various cohomologies we deal with. It provides cyclic bicomplexes computing our cohomologies, in which $S$ can be described as a shift.
As in cyclic cohomology, one introduce  these bicomplexes directly,
 and prove all the formulas in a straightforward manner. Here we prefer to apply the
formal constructions of \ref{Soper} to $W(\Ha)$ and to compute its $X$-complex. This computation can be carried out exactly as in the case of the tensor algebra (see Example \ref{tensor}), and this is done in the proof of Theorem \ref{calcul}. We end up with the following exact sequence of complexes (which can be taken as a definition):}
\[ .\ .\ .\ \rmap W^b(\Ha) \stackrel{t- 1}{\rmap} W(\Ha)\stackrel{N}{\rmap}W^b(\Ha) \stackrel{t- 1}{\rmap} W(\Ha)\rmap .\ .\ .\ ,\]
\emph{where we have to explain the new objects. First of all, $W^b(\Ha)$ is the same as $W(\Ha)$ but with a new boundary $b= \del + b_t$ with $b_t$ described below. The $t$ operator is the backward cyclic permutation:}
\[ t(ax)= (-1)^{|a||x|} xa ,\]
\emph{for $a\in H$ or of type $\omega_h$. This operator has finite order in each degree of $W(\Ha)$: 
we have $t^{p}= 1$ on elements of bi-degree $(p, q)$, so $t^{k!}= 1$ on elements of total 
degree $k$. The norm operator $N$ is $N:\, = 1+ t+ t^2+ . . .  +t^{p-1}$ on elements of 
bi-degree $(p, q)$.
The boundary $b$ of $W^b(\Ha)$ is $b= \del+ b_t$,}
\[ b_t(ax)= t(\del_0(a)x) ,\]
\emph{for $a\in \Ha $ or of type $\omega_h$. For all the operators
  involved, see also section \ref{grea}. Obviously, the powers $I(\Ha)^{n+ 1}$ are invariant by $b, t-1, N$, 
so we get similar sequences for $I_n(\Ha)$, $W_n(\Ha)$.\\
For reference, we conclude:}\end{num}

\begin{cor}\label{iii}There are exact sequences of complexes:\end{cor}
\begin{equation}\label{dac}
 CC(W_{n}(\Ha)):\ \ .\ .\ .\ \rmap W_{n}^b(\Ha) \stackrel{t- 1}{\rmap}
 W_{n}(\Ha)\stackrel{N}{\rmap}
W_{n}^b(\Ha) \stackrel{t- 1}{\rmap} 
W_{n}(\Ha)\rmap \ .\ .\ .\  
\end{equation}
\begin{equation}\label{seq1}
0 \rmap W_n(\Ha)_{\natural}
 \stackrel{N}{\rmap} W_{n}^b(\Ha) 
 \stackrel{t- 1}{\rmap} W_{n}(\Ha) 
 \stackrel{N}{\rmap} W_{n}^b(\Ha) 
 \stackrel{t- 1}{\rmap} W_{n}(\Ha) 
 \rmap \  .\ .\ . \end{equation}
\begin{equation}\label{seq2}
0\rmap I_n(\Ha)_{\natural}\  
 \stackrel{N}{\rmap} I_{n}^b(\Ha)\ \
 \stackrel{t- 1}{\rmap} I_n(\Ha)\ \
 \stackrel{N}{\rmap}I_{n}^b(\Ha)\ \
 \stackrel{t- 1}{\rmap} I_n(\Ha)\ 
 \rmap  \ .\ .\ .
\end{equation}

\begin{cor}\label{cordoi} There are short exact sequences of complexes:\end{cor}
\begin{equation}\label{seq3}
0\rmap W_n(\Ha)_{\natural} \stackrel{N}{\rmap}
W_{n}^{b}(\Ha)\stackrel{t- 1}{\rmap}
W_n(\Ha)\rmap
W_n(\Ha)_{\natural}\rmap 0 
\end{equation}
\begin{equation}\label{seq4}
0\rmap I_n(\Ha)_{\natural} \stackrel{N}{\rmap}\
I_{n}^{b}(\Ha)\stackrel{t- 1}{\rmap}\ \
I_n(\Ha)\rmap\ \
I_n(\Ha)_{\natural}\ \rmap 0 
\end{equation}\\

In particular, (\ref{seq1}), (\ref{seq2}), give bicomplexes which
compute the cohomologies of $W_n(\Ha)_{\natural}$,
$I_n(\Ha)_{\natural}$. They are similar to the 
(first quadrant) cyclic bicomplexes appearing in cyclic cohomology,
and come equipped with an obvious shift operator, which defines
our $S$- operation:
\[ S: H^{*}(W_n(\Ha)_{\natural})\rmap H^{*+ 2}(W_n(\Ha)_{\natural}),\]
(and similarly for $I_n(\Ha)_{\natural}$). Alternatively, one can 
obtain $S$ as cup-product by the $Ext^2$ classes arising from Corollary \ref{cordoi}.
 

\emph{End of proof of theorem \ref{th1}:} Denote for simplicity by $\underline{CC}^*(I_n), \underline{CC}^*(W), \underline{CC}^*(W_n)$ the (first quadrant) cyclic bicomplexes (or their total complexes) of $I_n$, $W$, and $W_n$, respectively. We have a map of short exact sequences of complexes:
\[ \xymatrix {
0 \ar[r] & \tilde{I}_{n}(\Ha)_{\natural} \ar[r] \ar[d]^-{N} & W(\Ha)_{\natural} \ar[r] \ar[d]^{N} & W_n(\Ha)_{\natural} \ar[r] \ar[d]^{N} & 0\\
0\ar[r] & \underline{CC}^*(I_n) \ar[r] & \underline{CC}^*(W) \ar[r] & \underline{CC}^*(W_n) \ar[r] & 0 } \]
where we have used the fact that $N: I_n(\Ha)_{\natural}\rmap I_n(\Ha)$ factors through the projection $I_n(\Ha)_{\natural}\twoheadrightarrow \tilde{I}_{n}(\Ha)_{\natural}$ (being defined on the entire $W(\Ha)_{\natural}$). Applying the five lemma to the exact sequences induced in cohomology by the previous two short exact sequences, the statement follows. \ \ $\Boxe$


\begin{num}{\bf Example:}\label{compus}\emph{ There are canonical
    Chern and Chern-Simons classes induced by any group-like element
    $\rho \in \Ha$ (i.e. with the property $\Delta(\rho)=
    \rho\otimes\rho$). Denote by $\omega$ its curvature. Since $\del 
    (\omega^n)=[\, \omega^n,\rho]$ is a commutator, $\omega^n$ 
define cohomology classes:}
\begin{equation} 
ch_{2n}(\rho):\ = [\, \natural ( \frac{1}{n!}\omega^n)] \in
H^{2n}(I_{m}(\Ha)_{\natural})\ ,
\end{equation}
\emph{for any $n\geq m$. The associated Chern-Simons class $cs_{2n-1}(\rho):\, = CS(ch_{2n}(\rho))$ is given by the formula (see (\ref{sim})):}
\[ cs_{2n-1}(\rho) = [\, \natural \{ \frac{1}{(n-1)!} \int_0^1 \rho (t \del(\rho)+ t^2\rho^{2})^{n-1} dt \}] \in H^{2n-1}(W_{m}(\Ha)_{\natural}) .\]
\emph{To compute $S(ch_{2n}(\rho))$, we have to solve successively the equations:}
\[ \left\{ \begin{array}{lll}
\del(\frac{1}{n!}\omega^n) & = & (t-1)(v) \\
b(v) & = & N(w)  \end{array}
\right. \]
\emph{and then $S(ch_{2n}(\rho))= [\, \natural(w)]$. The first equation has the obvious solution $v= \frac{1}{n!}\rho\omega^n$, whose $b(v)= \frac{1}{n!}\omega^{n+1}$, so the second equation has the solution $w= \frac{1}{(n+1)!}\omega^{n+1}$. In conclusion,}
\begin{equation}\label{Srel}
S(ch_{2n}(\rho))= ch_{2\ps n+1\pd}(\rho), \ \ S(cs_{2n-1}(\rho)) = cs_{2n+1}(\rho) .\end{equation}
\emph{(where the second relation follows from the first one and Theorem \ref{th1}.)}
\end{num}

\section{The Weil complex and higher traces}
\label{NCCW}

\hspace*{.3in}We explain now how the Weil complex introduced in the previous section appears naturally in the case of higher traces, and Hopf algebra actions. The main reason that $HC_{\delta}^*(\Ha)$ is still the target of these characteristic maps is that it can be computed by the truncation of the Weil complex (see Theorem \ref{comp}, whose proof is postponed until the next section). To prove the compatibility with the $S$-operator, we first have interpret the complexes introduced in the previous section in terms of Cuntz-Quillen's (tower of) relative $X$-complexes. We will obtain in particular the case of usual traces discussed in Section \ref{Hocic}. Also, for $\Ha= {\bf C} \rho$ (example \ref{simple}), we re-obtain the results , and interpretations of some of the computations of \cite{Qext} (see Example \ref{eli} below).\\
\hspace*{.3in}In this section $\Ha$ is a Hopf algebra, $\delta$ is a character on
$\Ha$, and $A$ is a $\Ha$-algebra. We assume for simplicity that $S_{\delta}^2= Id$. \\

\begin{num}\label{equiv}{\bf Localizing $W(\Ha)$:} \emph{First of all remark that the 
Weil complex $W(\Ha)$ is naturally an $\Ha$ DG algebra. By this we mean a DG algebra, 
endowed with a (flat) action, compatible with the grading and with the differentials. 
The action is defined on generators by:}
\[ g\cdot i\,(h):\ = i\,(gh),\ \ g\cdot \omega_h:\ = \omega_{gh},\ \ \forall\  g, h\, \in  \Ha .\]
\emph{and extended by $h(xy)=\sum h_0(x)h_1(y)$. Here,  to avoid confusions, 
we have denoted by $i:\Ha \rmap W(\Ha)$ the inclusion. Remark that the action 
preserves the bi-degree (see \ref{bigrad}), so $W(\Ha)_{\delta}$ has an induced 
bi-grading.  We briefly explain how to get the localized version for the constructions 
and the properties of the previous section. First of all one can localize with respect 
to $\delta$ as in Section $3$, and (with the same proof as of Proposition \ref{descend}), 
all the operators descend to the localized spaces. The notation $I_n(\Ha)_{\natural, \delta}$ 
stands for $I_n(\Ha)$ divided out by commutators and co-invariants. For Theorem \ref{th1}, 
remark that the contraction used there is compatible with the action. 
To get the exact sequences from Corollary \ref{iii} and \ref{cordoi}, we may look at 
them as a property for the cohomology of finite cyclic groups, acting (on each fixed bi-degree)
 in our spaces. Or we can use the explicit map $\alpha: W(\Ha)\rightarrow W(\Ha)$ defined by 
$\alpha := (t+ 2t^2+ . . . + \ps p-1\pd t^{p-1})$ on elements of bi-degree $(p, q)$, which has
 the properties: $(t-1)\alpha+ N= pId$, $\alpha(I(\Ha)^{n+1})\subseteq I(\Ha)^{n+1}$, 
and $\alpha$ descends (because $t$ does). So, also the analogue of Theorem \ref{th1} follows. 
In particular $H^{*}_{\delta}(W_n(\Ha)_{\natural})$ is computed either by the complex 
$W_n(\Ha)_{\natural, \delta}$, or by the (localized) cyclic bicomplex 
$CC_{\delta}^*(W_n(\Ha))$ (analogous to (\ref{dac})). Similarly, we consider the $S$ operator, and the periodic versions of 
these cohomologies. Due to the shift in the degree already existent in
the case of the tensor algebra (see \ref{quillen}), we re-index these
cohomologies:}
\end{num}

\begin{defin} Define $HC_{\delta}^*(\Ha, n):\ = \
  H^{*+1}(W_n(\Ha)_{\natural, \delta})$, and denote by $CC^{*}_{\delta}(\Ha,
  n)$ the cyclic bicomplex computing
it, that is, $CC_{\delta}^*(W_n(\Ha))$ shifted by one in the vertical
  direction.
\end{defin}

Remark that for $n=0$ we obtain Connes-Moscovici's cyclic cohomology and:
\[ CC^*(\Ha, 0) = CC^*(\Ha)\ , \ \ CC_{\delta}^*(\Ha, 0) = CC_{\delta}^*(\Ha), \]
 while, in general, there are obvious maps:
\begin{equation}\label{tower} .\ .\ .\ \stackrel{\pi_3}{\rmap} HC_{\delta}^*(\Ha, 2)
\stackrel{\pi_2}{\rmap} HC_{\delta}^*(\Ha, 1) 
\stackrel{\pi_1}{\rmap} HC_{\delta}^*(\Ha, 0) \cong HC_{\delta}^*(\Ha) .
\end{equation}
In the next section we will prove:
\begin{st}\label{comp} $HC_{\delta}^{*}(\Ha, n)\cong
  HC_{\delta}^{*-2n}(\Ha)$, and the tower (\ref{tower}) is the
$S$ operation tower for $HC_{\delta}^{*}(\Ha)$. More precisely, there
are isomorphisms $\beta: HC_{\delta}^{*}(\Ha, n)\tilde{\rmap}
HC_{\delta}^{*-2}(\Ha, n-1)$
such that the following diagram commutes:
\[ \xymatrix{
.\ .\ .\ \ar[dr]^-{\pi} & 
.\ .\ .\ \ar[dr]^-{\pi} \ar[d]^-{\bf{\beta}} &
.\ .\ .\ \ar[dr]^-{\pi} \ar[d]^-{{\bf \beta}}&
.\ .\ .\ \ar[dr]^-{\pi} \ar[d]^-{\beta} & 
.\ .\ . \\
.\ .\ .\ \ar[dr]^-{\pi}  \ar[r] &
HC_{\delta}^{*}(\Ha,2)\ar[d]^-{\beta}\ar[r]^-{S} \ar[dr]^-{\pi}    &
HC_{\delta}^{*+2}(\Ha,2)\ar[d]^-{\beta} \ar[r]^-{S}\ar[dr]^-{\pi}    &
HC_{\delta}^{*+4}(\Ha,2)\ar[d]^-{\beta} \ar[r]^-{S}\ar[dr]^-{\pi}    &
.\ .\ . \\
.\ .\ .\ \ar[dr]^-{\pi} \ar[r] &
HC_{\delta}^{*-2}(\Ha,1)\ar[d]^-{\beta}\ar[r]^-{S} \ar[dr]^-{\pi}    &
HC_{\delta}^{*}(\Ha,1)\ar[d]^-{\beta} \ar[r]^-{S}\ar[dr]^-{\pi}    &
HC_{\delta}^{*+2}(\Ha,1)\ar[d]^-{\beta} \ar[r]^-{S}\ar[dr]^-{\pi}    &
.\ .\ . \\
.\ .\ .\    \ar[r] &
HC_{\delta}^{*-4}(\Ha,0)\ar[r]^-{S}     &
HC_{\delta}^{*-2}(\Ha,0) \ar[r]^-{S}   &
HC_{\delta}^{*}(\Ha,0)\ar[r]^-{S}  &
.\ .\ . } \]
\end{st}


\begin{num}\label{daev}{\bf The case of even equivariant traces:}
  \emph{ Consider now an equivariant even trace over $A$,  i.e. an extension:
\begin{equation}\label{exten} 0\rmap I\rmap R\stackrel{u}{\rmap} A\rmap 0 \end{equation}
 and a $\delta$-invariant trace $\tau : R\rmap \bf{C}$ vanishing on $I^{n+1}$. 
To describe the induced characteristic map, we choose a linear splitting $\rho:
A\rmap R$ of (\ref{exten}). As in the case of the usual Weil complex, 
there is a unique equivariant map of DG algebras:}
\[ \tilde{k} : W(\Ha) \rmap Hom(B(A), R), \]
\emph{sending $1\in \Ha$ to $\rho$. This follows from Proposition \ref{univ} and from 
the equivariance condition (with the same arguments as in
\ref{caract}). Here, the action of $\Ha$ on $Hom(B(A), R)$ is induced by 
the action on $R$. Since $\rho$ is a homomorphism modulo $I$,
$\tilde{k}$ sends $I(\Ha)$ to $Hom(B(A), I)$, so induces a map
$W_n(\Ha)\rmap Hom(B(A), R/I^{n+ 1})$. As in \ref{caract}, composing with the 
$\delta$-invariant trace:}
\[ \tau_{\natural}: Hom(B(A), R/I^{n+1})\rmap C_{\lambda}^*(A)[1], \ \ \phi \mapsto  \tau\compose \phi\compose N\ , \]
\emph{we get a $\delta$-invariant trace on $W_n(\Ha)$, so also a chain map:}
\begin{equation}
\label{care} k^{\tau, \rho}: W_n(\Ha)_{\natural, \delta} \rmap
C_{\lambda}^*(A)[1]\ .
\end{equation}
\emph{Denote by the same letter the map induced in cohomology:}
\begin{equation}\label{charact}  
k^{\tau, \rho}: HC_{\delta}^*(\Ha, n)\rmap HC^{*}(A), 
\end{equation}
\emph{or, using the isomorphism of Theorem \ref{comp}:}
\begin{equation}
\label{characti}  k^{\tau, \rho}: HC_{\delta}^{*- 2n}(\Ha)\rmap HC^{*}(A), 
\end{equation}
\end{num}

\begin{st}\label{damm} The characteristic map (\ref{characti}) of the even higher trace $\tau$ does not depend on the choice of the splitting $\rho$ and is compatible with the $S$-operator.\end{st}

\emph{proof:}(compare to \cite{Qext}) We use \ref{path}. If $\rho_0, \rho_1$ are two liftings, form $\rho= t\rho_0+ (1-t)\rho_1 \in Hom(A, R[t])$, viewed in the degree one part of the DG algebra $Hom(B(A), R\otimes\Omega(1))$. It induces a unique map of $\Ha$DG algebras $\tilde{k}_{\rho}: W(\Ha)\rmap Hom(B(A), R\otimes\Omega(1))$, sending $1$ to $\rho$, which maps 
$I(\Ha)$ to the DG ideal $Hom(B(A), I\otimes\Omega(1))$ (since $\omega_{\rho}$  belongs to the former). Using the trace $\tau\otimes\int : R/I^{n+1}\otimes\Omega(1)\rmap \bf{C}$, and the universal cotrace on $B(A)$, it induces a chain map:
\[ k_{\rho_0, \rho_1}: W_n(\Ha) \rmap C_{\lambda}^*(A)[1],\]
which kills the coinvariants and the commutators. The induced map on 
$W_n(\Ha)_{\natural, \delta}$ is a homotopy between $k^{\tau, \rho_0}$ and $k^{\tau, \rho_1}$. 
The compatibility with $S$ follows from the fact that the characteristic map (\ref{care}) 
can be extended to a map between the cyclic bicomplexes $CC_{\delta}^*(\Ha, n)$ and $CC^*(A)$. We will prove this after shortly discussing the case of odd higher traces.
\ \ $\Boxe$ 

\begin{num}\label{daod}{\bf The case of odd equivariant traces:}
\emph{ A similar discussion applies to the case of odd equivariant 
traces on $A$, i.e. extensions 
(\ref{exten}) endowed with a $\delta$-invariant linear map 
$\tau : I^{n+1}\rmap \bf{C}$,
vanishing on $[I^n, I]$ . The
resulting map $H^{*+1}(I_n(\Ha)_{\natural, \delta})\rmap HC^{*-1}(A)$,
combined with Corollary \ref{th1}, the comments in \ref{equiv},
and Theorem \ref{comp}, give the characteristic map:}
\begin{equation}\label{zbir} k^{\tau, \rho} : HC_{\delta}^{*- 2n- 1}(\Ha)\cong
HC^{*-1}_{\delta}(\Ha, n) \rmap HC^{*}(A),\end{equation}
\emph{which has the same properties as in the even case:} \end{num}

\begin{st}\label{dam} The characteristic map (\ref{zbir}) of the odd
  higher trace $\tau$ does not depend on the choice of the splitting
  $\rho$ and is compatible with the $S$-operator.\end{st}

\begin{num}{\bf The localized tower $\xx_{\delta}(R, I)$:}\emph{
    Recall that given an ideal $I$ in the algebra $R$, one has a tower
    of super-complexes $\xx _{\delta}(R, I)$ given by  (\cite{CQ3}, pp. 396):}
\[ \xymatrix{ \xx^{2n+ 1}(R, I): R/I^{n+1} \ar@<-1ex>[r]_-{d} & \Omega^1(R)_{\natural}/\natural(I^{n+ 1}dR+ I^n dI) \ar@<-1ex>[l]_-{b}}\ ,\]
\[ \xymatrix{ \xx^{2n}(R, I): R/(I^{n+1}+ [I^n, R]) \ar@<-1ex>[r]_-{d} & \Omega^1(R)_{\natural}/\natural(I^{n} dR)\ar@<-1ex>[l]_-{b}} ,\]
\emph{where $\natural : \Omega^1(R)\rightarrow \Omega^1(R)_{\natural}$ is the projection.  The structure maps $\xx^{n}(R, I)\rightarrow \xx^{n+1}(R, I)$ of the tower are the obvious projections. We have a localized version of this, denoted by $\xx_{\delta}(R, I)$, and which is defined by:}
\[ \xymatrix{ \xx^{2n+ 1}_{\delta}(R, I): R/(I^{n+1}+ {\rm coinv}) \ar@<-1ex>[r]_-{d} & \Omega^1(R)_{\natural, \delta}/\natural(I^{n+ 1}dR+ I^n dI) \ar@<-1ex>[l]_-{b}}\ ,\]
\[ \xymatrix{ \xx^{2n}_{\delta}(R, I): R/(I^{n+1}+ [I^n, R]+ {\rm coinv}) \ar@<-1ex>[r]_-{d} & \Omega^1(R)_{\natural, \delta}/\natural(I^{n} dR)\ar@<-1ex>[l]_-{b}} ,\]
\emph{where, this time,  $\natural$ denotes the projection $\Omega^1(R)\rightarrow \Omega^1(R)_{\natural, \delta}$.\\
\hspace*{.3in} Remark that the construction extends to the graded
case, and each $\xx^n(R, I)$ is a super-complex of complexes.}\end{num}

\begin{st}\label{calcul} The cyclic bicomplex $CC_{\delta}^*(\Ha, n)$ is isomorphic to the 
bicomplex \linebreak $\xx_{\delta}^{2n+1}(W(\Ha), I(\Ha))$. \end{st}

\emph{proof:} The computation is similar to the one of $X(T\Ha)$ (see Example \ref{tensor} and 
Proposition \ref{bar}). Denote $W= W(\Ha)$, $I= I(\Ha)$, and let $V\subset W(\Ha)$ be the linear 
subspace spanned by $h$'s and $\omega_h$'s. Remark that $W$, as a graded algebra, is freely 
generated by $V$. This allows us to use exactly the same arguments as in \ref{tensor}, \ref{bar}
to conclude that $\Omega^1W\cong \tilde{W}\otimes V\otimes\tilde{W},$\ $\Omega^1(W)_{\natural}
\cong V\otimes \tilde{W}= W,$\ $\Omega^1(W)_{\natural, \delta}\cong W_{\delta}$; the 
projection $\natural : \Omega^1(W)\rightarrow \Omega^1(W)_{\natural}$ identifies with:
\begin{equation}\label{proj} \natural : \Omega^1(W)\rmap V\otimes\tilde{W} = W,\ \  
x \partial_u(v)y\mapsto (-1)^{\mu} vyx\ ,\end{equation}
for $x, y\in\ \tilde{W}, v\in V$. Here $\mu= deg\ps x\pd(deg\ps v\pd + deg\ps y\pd)$ introduces 
a sign, due to our graded setting, and $\partial_u: W\rightarrow \Omega^1(W)$ stands for
the universal derivation of $W$. Using this, we can compute the new boundary of $W$, coming
from the isomorphism $W\cong \Omega^1(W)_{\natural}$, and we end up we the $b$-boundary of $W$,
defined in Section $5$. For instance, if $x= hx_0\in W$ with $h\in \Ha$, since 
$\natural (\partial_u(h)x_0)= x$ by (\ref{proj}), its boundary is:
\begin{eqnarray*}
\lefteqn{\natural(\partial_u(\del(h)x_0)- \partial_u(h)\del(x_0)) =}\\
 & & = \natural( \partial_u (\omega_h- \sum h_0h_1)x_0- \partial_u(h) \del(x_0)) = \\
 & & = \natural( \partial_u(\omega_h)x_0- \sum \partial_u(h_0)h_1x_0 - \sum h_0\partial_u(h_1)x_0 - 
           \partial_u(h) \del(x_0))= \\
 & & = \omega_hx_0 - \sum h_0h_1x_0 - \sum (-1)^{deg\ps x\pd} h_1x_0h_0 - h\del(x_0) = \\
 & & = \del(h)x_0 - \sum t(h_0h_1x_0) - h\del(x_0) =  \\
 & & = \del(hx_0) +
 t(\del_0(h)x_0)= b(hx_0) 
\end{eqnarray*}
Remark also that our map (\ref{proj}) has the property:
\begin{equation}\label{ideal} \natural(I^n\partial_u I + I^{n+1}\partial_u W)= I^{n+1}.
\end{equation}
These give the identification $\xx^{2n+1}(W, I)\cong CC^*(\Ha, n)$. The localized 
version of this is just a matter of checking that the isomorphism 
$\Omega^1(W)_{\natural, \delta}
\cong W_{\delta}$ already mentioned, induces $\Omega^1(W)_{\natural, \delta}
/ \natural(I^n \partial_u I + I^{n+1}\partial_u W) \cong (W/I^{n+1})_{\delta}$, which follows from
(\ref{ideal}). \ \ $\Boxe$ \\

\begin{num} {\bf Proof of the $S$-relation:}\emph{ We freely use the
    dual constructions for (DG) coalgebras $B$, 
such as the universal coderivation $\Omega_1(B)\rightarrow B$, the space of co-commutators
$B^{\natural}= Ker(\Delta - \sigma\compose\Delta: B\rightarrow B\otimes B)$, and the $X$-complex
$X(B)$ (see \cite{Qext}). Denote $B= B(A)$,
$L =Hom(B, R)$, $J= Hom(B, I)$. Our goal is to prove that the characteristic map (\ref{care}) can
 be defined at
the level of the cyclic bicomplexes. Consider first the case of even traces $\tau$.
Since the $\Ha$ DG algebra map 
$\tilde{k}: W(\Ha)\rightarrow L$ maps $I(\Ha)$ 
inside $J$, there is an induced  map $\xx_{\delta}^{2n+1}(W(\Ha), I(\Ha))
\rightarrow \xx_{\delta}^{2n+1}(L, J)$, extending 
$W_n(\Ha)_{\natural, \delta}\rightarrow (L/J^{n+1})_{\natural, \delta}$. So, it suffices to
show that the map $(L/J^{n+1})_{\natural, \delta}\rightarrow Hom(B^{\natural}, (R/I^{n+1})_{\natural, \delta})$
(constructed as (\ref{care})), lifts to a map of super-complexes (of complexes)}
\begin{equation}\label{exte} 
\xx_{\delta}^{2n+1}(L, J)\rightarrow Hom(X(B), (R/I^{n+1})_{\natural, \delta})
\end{equation}
\emph{Indeed, using Theorem \ref{calcul}, the (similar) computation of $X(B)$ (as the cyclic 
bicomplex of $A$), the interpretation of the norm map $N$ as the universal cotrace 
of $B$ (see \cite{Qext}), and the fact that any $\tau$ 
as above factors
through $(R/I^{n+1})_{\natural, \delta}\rightarrow {\bf C}$, the map (\ref{exte})
is 'universal' for our problem.
 The construction of  (\ref{exte}) is quite simple. The composition
with the universal coderivation of $B$ is a derivation
$L\rightarrow Hom(\Omega_1(B), R)$ on $L$, so it induces a map $\chi: \Omega^1(L)\rightarrow
Hom(\Omega_1(B), R)$. Since $\chi$ is a $L$-bimodule map, and it is compatible with the action of 
$\Ha$, it induces a map $\Omega^1(L)_{\natural}\rightarrow Hom(\Omega_1(B)^{\natural},
(R/I^{n+1})_{\natural, \delta})$, which kills $\natural(J^{n}dJ+ J^{n+1}dL+ coinv)$. This,
together with  
with the obvious $(L/J^{n+1})_{\delta}\rightarrow Hom(B, (R/I^{n+1})_{\natural, \delta})$, 
give (\ref{exte}). For the case of odd higher traces we proceed
similarly, and use the remark that
(\ref{exte}) was apriori defined at the level of $L$, and $\Omega^1(L)$, and
one can restrict to the ideals (instead of dividing out by them).\ \ $\Boxe$ }\end{num}  

\begin{num}\label{eli}{\bf Examples:}\emph{ Choosing $\rho= 1\in \Ha$
    (the unit of $\Ha$) in Example \ref{compus}, and applying the
    characteristic map to the resulting classes, we get the
    Chern/Chern-Simons classes (in the cyclic cohomology of $A$),
    described in \cite{Qext}. Remark that our proof of the
    compatibility with the $S$ operator consists on two steps: the
    first one proves the universal formulas (\ref{Srel}) at the level
    of the Weil complex, while the second one shows, in a formal way,
    that the characteristic map can be defined at the level of the
    cyclic bicomplexes. This allows us to avoid the explicit cochain
    computations.\\
Another interesting example is when $\Ha= U(\lie)$ as in Example
    \ref{exemple}, $\delta=$ the
    counit. Via the computation of Theorem \ref{env}, our construction
    associates to any $G$-algebra $A$, and any $G$-invariant
    higher trace $\tau$ on $A$, of parity $i$, a $\bf{Z}/2\bf{Z}$ graded characteristic maps:}
\begin{equation}\label{exca}
k_{\tau}: H_*(\lie)\rmap HP^{*+i}(A).
\end{equation}  
\emph{When $\tau$  is a usual invariant trace $\tau: A\rmap {\bf C}$,
  we have the following formula (use (\ref{reftau}) and the map $A$
  used in the proof of Theorem \ref{env}):
\begin{eqnarray*}
\lefteqn{k_{\tau}(v_1\wedge .\, .\, . \wedge v_n)(a_0, a_1, .\, .\, .\, , a_n)=}\\
& & = \frac{1}{n!}\sum_{\sigma} sign(\sigma) \tau( a_0 v_{\sigma\ps 1\pd}(
a_1) \, .\, .\, .\,  v_{\sigma\ps n\pd}(a_n))\ ,
\end{eqnarray*} 
( where $v(a):= L_v(a)$ is the Lie derivative). This coincides 
with the characteristic map described in \cite{Coo}.}\end{num}

\section{Proof of Theorem \ref{comp}, and equivariant cycles}
\label{grea}

\hspace*{.3in}This section is devoted to the proof of Theorem \ref{comp}. At the end we illustrate how the new complexes computing $HC_{\delta}^*(\Ha)$ which arise during the proof can be used to construct characteristic maps associated to equivariant cycles.\\
\hspace*{.3in}We first concentrate on the non-localized version, whose proof
uses  explicit formulas which can be  easily 
 localized. So,  we construct
isomorphisms 
\[\beta: H^*(W_n(\Ha)_{\natural})\tilde{\rmap}H^{*-2}(W_{n-1}(\Ha)_{\natural})\]
(and explicit inverses)  such that the following diagram is
commutative:
\[ \xymatrix{ 
H^*(W_{n, \natural}) \ar[r]^-{S} \ar[d]_-{\beta}\ar[rd]^-{\pi} &
H^{*+ 2}(W_{n, \natural})\ar[d]^-{\beta}\\
H^{*- 2}(W_{n- 1, \natural})\ar[r]^-{S} &
H^*(W_{n- 1, \natural})}\]

Let's start by fixing some notations. Denote $W_n(\Ha)= W_n$, $I(\Ha)= I$, $I^{\ps n\pd}=
I^{n+1}/I^n$ viewed as the subspace of $W$ spanned by elements having
exactly $n$ curvatures. The only grading we consider is the total
grading (with $deg(h)=1, deg(\omega_h)= 2$); notations like
$(W_n,\del), (W_n, b)$ are used to specify the complexes we are
working with. In general, if the (signed) cyclic permutation acts on
a vector space $X$, denote $X_{\natural}= X/Im(1-t)$.\\
We review now the various operators we have. First of all, 
\[ \del= \del_0+ d,\]
where $\del_0, d$ are the degree $1$ derivations given on generators
by:
\[ \del_0(h)= -\sum h_0h_1, \ \  \del_0(\omega_h)= -\sum \omega_{h_0}h_1-
h_0\omega_{h_1},\]
\[ d(h)= \omega_h,\ \ d(\omega_h)= 0. \]
Secondly, the operator $b= \del+ b_t$, where:
\[ b_t(h x)= \sum (-1)^{deg x} h_1xh_0,\ \  b_t(\omega_h x)=\sum
h_1x\omega_{h_0}- (-1)^{deg x} \omega_{h_1}xh_0.\]
Define also $b_0= \del_0+ b_t$. It is straightforward to check: 
\[ b= b_0+ d,\ \  b_0^2= d^2= [b_0, d]= 0.\]
Point out that $d$ commutes with $t$.\\
For the construction of $\beta$ we need the following degree $-1$ operator:
\[\theta: W\rmap W,\ \ \theta(hx)=0,\ \ \theta(\omega_hx)=hx .\]
For constructing the inverse of $\beta$, we will use the degree $0$
operators $\phi_i:I\rmap W$, $i=0, 1$. On $I^{\ps n\pd}$,
\[ \phi_1(hx)= 0, \ \ n \phi_1(\omega_h x)= x\omega_h . \] 
For $y= x_0x_1 .\, .\, .\, x_p \in I^{\ps n\pd}$, where each of the
$x_i$'s are of type $h$ or $\omega_h$, we put $\lambda_i(y)=\#\{ j\leq i:
x_j {\rm\ is\ of\ type\ } \omega_h\}$, and define  $n \phi_0(y)=
\sum_{1}^{n-1} \lambda_i(y) t^{i}(y).$
For a conceptual motivation, see the next proof.
We can actually forget about these formulas, and just keep their
relevant properties:

\begin{lem}\label{formulas} (i) $[\theta, b_0]= 0,\ [\theta, \del_0]=
  0,\ [\theta, d]=1,\ \theta^2= 0 $,\\
(ii) $\phi_1N- (1-t)\phi_0= 1, N\phi_1- \phi_0(1-t)= 1$,\\  
(iii) $\phi_1b_0= \del_0\phi_1$ modulo $Im(1- t)$, $\phi_1 \theta= 0$, \\
(iv)  $\phi_0\del_0= b_0\phi_0$ modulo $Im \theta$.
\end{lem}

\emph{Proof:} (i) and (iii) follow by direct computation. For
instance, for the first part of (iii) one has $\phi_1b_0= 0= 
\del_0\phi_1$ on elements of type $hx$, while on elements of type
$\omega_h x\in I^{\ps n\pd}$:
\begin{eqnarray*}
\phi_1b_0 (\omega_h x)= & \phi_1((\omega_{h_0}h_1-
      h_0\omega_{h_1})x+ \omega_h\del_0(x)+ h_1x\omega_{h_0}-
      (-1)^{deg\ps x\pd}\omega_{h_1}xh_0)=\\
& = n( h_1 x\omega_{h_0} + \del_0 (x) \omega_h -(-1)^{deg\ps x\pd}xh_0 \omega_{h_1})
\end{eqnarray*}
\[ \del_0\phi_1(\omega_hx)= n\del_0(x\omega_h)= n( \del_0(x)\omega_h+
      (-1)^{deg\ps x\pd} x\omega_{h_0}h_1- (-1)^{deg\ps x\pd}xh_0
      \omega_{h_1}), \]
and the two expressions are clearly the same modulo $Im(1- t)$.\\ 
One can check  directly also (ii). Instead, let's 
explain that $\phi_0$, $\phi_1$ have been constructed in such a
way that (ii) holds. On the graded algebra $W= \oplus I^{\ps n\pd}$,
we have a Goodwillie \cite{Good}  type derivation: multiplication by the number 
of curvatures. Since $W$ is a tensor algebra, it comes equipped with a
canonical connection (see the end of section $3$ in \cite{CQ2}). We
know that the $X$- complex of $W$ is the cyclic bicomplex, and the
general Cartan homotopy formula of \cite{CQ3} for our derivation $D$,
gives precisely the homotopy $(n\phi_0, n\phi_1)$ on $I^{\ps  n\pd}$.
For (iv), remark first that $Im \theta= Ker \theta$, so it suffices to
show that $A:=(\theta \phi_0)\del_0- b_0(\theta \phi_0)$ is zero. From the
first formula of (iii), the second of (ii), and (i), it follows
that $A (1- t)=0$. So, it is enough to check $A= 0$ on homogeneous monomials
having a curvature as first element. Such an element can be written 
as $X= \omega(h^1)X^1\, .\, .\, . \, \omega({h^n})X^n \in I^{\ps n\pd}$, where $X^i\in
\widetilde{T}(\Ha)$, $\omega(h)= \omega_h$. On $X$, $\theta \phi_0(X)= \sum_{1}^n
\epsilon_i (i-1) h^iX^i\omega(h^{i+1})X^{i+1}\, .\, .\, . \, \omega(h^n)X^n
\omega(h^1)X^1\, .\, .\, . \,\omega(h^{i-1})X^{i-1}$ ($\epsilon_i$ are
corresponding signs), and the
proof becomes a lengthy straightforward computation.\ \ \ $\Boxe$\\

To define $\beta$, we need the right complexes computing $H^*(W_{n,
  \natural})$. One of them is given by the following:

\begin{lem}\label{ncom} (i) There are isomorphisms $p: H^*(W_{n, \natural})
  \tilde{\rmap}  H^*(I^{\ps n\pd}_{\natural}/Im d, \del_0)$,
  compatible with the $S$-operations.\\
\hspace*{.3in} (ii) One has short exact sequences:
\begin{equation}\label{nustiu}
0\rmap (I^{\ps n-1\pd}_{\natural}/Im d, \del_0) \stackrel{d}{\rmap} (I^{\ps
  n\pd}_{\natural}, \del_0) \rmap (I^{\ps n\pd}_{\natural}/Im d,
  \del_0)[1]\rmap 0,
\end{equation}
whose induced boundary in cohomology identifies, via $p$, with $\pi$:
\[ \xymatrix{
H^*(W_{n, \natural}) \ar[r]^-{\pi} \ar[d]_-{p} & 
H^*(W_{n-1,\natural})\ar[d]_-{p}\\
H^*(I^{\ps n\pd}_{\natural}/Im d, \del_0)\ar[r]^-{\delta}&
 H^* (I^{\ps  n-1\pd}_{\natural}/Im d, \del_0)} \]
\end{lem}

Here, we view $(W_{n, \natural}, \del)$ as the total complex of the
double complex:

\[ 0\rmap (I^{\ps 0\pd}_{\natural}, \del_0) \stackrel{d}{\rmap}(I^{\ps 1\pd}_{\natural}, \del_0) \stackrel{d}{\rmap}
 .\ .\ .\ \stackrel{d}{\rmap} (I^{\ps n\pd}_{\natural}, \del_0) \rmap
 0\rmap .\, .\, .\]
and  $p$ is induced by the obvious
 augmentation sending $[\sum_0^n x_i]$ into $[x_n]$ ($x_i\in I^{\ps
 i\pd}$). The $S$ operation on $H^*(I^{\ps n\pd}_{\natural}/Im d,
 \del_0)$ of (i) is defined by the  cyclic bicomplex
 which is augmentation of (\ref{seq1}),
or, similar to (\ref{seq3}), by the $Ext^2$ class defined by the extension:
\begin{equation}\label{exten1}
0\rmap (I^{\ps n\pd}_{\natural}/Im d, \del_0) \stackrel{N}{\rmap}
(I^{\ps n\pd}/Im d, b_0) \stackrel{t- 1}{\rmap}
(I^{\ps n\pd}/Im d, \del_0) \rmap
(I^{\ps n\pd}_{\natural}/Im d, \del_0) \rmap 0 
\end{equation}
Using that $W$ is contractible along $d$ (cf. \ref{formulas} (i)), 
(i) is clear. Using that $d\compose t= t\compose d$, and that  taking invariants under the action of
a finite group does not affect exactness, also the first part of (ii) follows, while the last part is a
routine spelling out of the boundary of long exact sequences.\\
There is a slight modification of (\ref{exten1}) which can be used
to compute $H^*(W_{n, \natural})$, obtained as follows: 
(\ref{exten1}) splits into two short exact sequence:
\begin{equation}\label{exten2}
0\rmap (I^{\ps n\pd}_{\natural}/Im d, \del_0) \stackrel{N}{\rmap}
(I^{\ps n\pd}/Im d, b_0) \rmap (I^{\ps n\pd}/Im d+ Im N, b_0) \rmap 0,
\end{equation}
\begin{equation}\label{exten3}
0\rmap 
(I^{\ps n\pd}/Im d+ Im N , b_0) \stackrel{t- 1}{\rmap}
(I^{\ps n\pd}/Im d, \del_0) \rmap
(I^{\ps n\pd}_{\natural}/Im d, \del_0) \rmap 0 
\end{equation}
Since the middle complex of (\ref{exten3}) is acyclic, e.g. by
using the contraction $s_{-1}(hx)= \epsilon(h)x$, $s_{-1}(\omega_hx)=
0$ (which commutes with $d$), we get a quasi-isomorphism (which, in cohomology, is independent
of the contraction):
\begin{equation}\label{exten4}
 s_{-1}(1- t): (I^{\ps n\pd}/Im d+ Im N , b_0) \stackrel{q. i.}{\rmap}
 (I^{\ps n\pd}_{\natural}/Im d, \del_0)[1]
\end{equation}
Via this, the $S$ operator is simply the boundary of the long exact
sequence induced by (\ref{exten2}). \\
Now, our map is defined as the chain map:
\begin{equation}\label{exten5}
\beta: (I^{\ps n\pd}_{\natural}/Im d, \del_0) \rmap (I^{\ps n- 1\pd}/Im d+ Im N , b_0))[1]
\end{equation}
induced by $-\theta\compose N$. To understand our choice of complexes,
let's just mention that (\ref{exten5}) is an isomorphism when $n=
1$. Note also that $\beta$, as well as the chain map $\alpha$ bellow
describing its homotopical inverse, do not depend on the structure of
$\Ha$, other then the vector space structure.\\
  Now, to see that $\beta$ is compatible with the $S$ operation, and 
to construct its inverse (in cohomology), we make use of  the
following diagram with exact rows and columns:

\[ \xymatrix{
(I^{\ps n- 1\pd}_{\natural}/Im d, \del_0)[1] \ar[r]^-{d}\ar[d]_-{N} &
(I^{\ps n\pd}_{\natural}, \del_0) \ar[r]^-{p_5} \ar@<-1ex>[d]_-{\tilde{N}} \ar[dl]_-{-\tilde{\beta}} &
(I^{\ps n\pd}_{\natural}/Im d, \del_0) \ar[d]^-{N} \\
(I^{\ps n- 1\pd}/Im d, b_0)[1]\ar@<+1ex>[r]^-{\tilde{d}}\ar[d] &
(I^{\ps n\pd}, b_0) \ar@<+1ex>[l]^-{\tilde{\theta}}
\ar@<-1ex>[u]_-{\tilde{\phi}_1}\ar@<+1ex>[r]^-{p_1}
\ar@<-1ex>[d]_-{p_2} & 
(I^{\ps n\pd}/Im d, b_0) \ar@<+1ex>[l]^-{r}\ar[d]^-{p_3}\ar[dl]^-{p_2\compose r} \\
(I^{\ps n- 1\pd}/Im d+ Im N , b_0)[1] \ar[r]_-{d}&
(I^{\ps n\pd}/Im N , b_0) \ar@<-1ex>[u]_-{s} \ar[r]_-{p_4} &
(I^{\ps n\pd}/Im d+ Im N , b_0)} \]

Here $\tilde{N}, \tilde{\phi}_1, \tilde{d}, \tilde{\theta}$ are induced
by $N, \phi_1, d, \theta$, respectively, $p_1, p_2$ are the obvious
projections, $r$ is the map induced by $\theta d$, $s$ is the
one induced by $\phi_0 (1- t)$, and $\tilde{\beta}$ is the one induced by
$-\theta N$. From
\ref{formulas} (ii), (iii), $\tilde{\phi}_1, s$ are chain maps with:
\begin{equation} \label{uli}
\tilde{\phi}_1 N= Id,\  p_2 s= Id,\  N \tilde{\phi}_1+ s p_2= Id .
\end{equation}
Also, from (i) of the same Lemma, $\tilde{\theta}, r$ are chain maps
with:
\begin{equation} \label{sumunu}
 \tilde{\theta} \tilde{d}= Id,\ p_1 r= Id, \ r p_1+ \tilde{d}
 \tilde{\theta}= Id. 
\end{equation}
Since $-\tilde{\beta} d= \tilde{\theta} N d= \tilde{\theta} \tilde{d}
N= N$, $\tilde{\beta}$ induces a map between the Cokernels of $d$ and
$N$, and this is precisely our map (\ref{exten5}). Moreover,
$\tilde{\beta}$ induces a map between the left vertical short exact
sequence, and the upper horizontal one. The boundaries of the long
exact sequences induced in cohomology are, by the previous remarks and
by (ii) of Lemma \ref{ncom}, the $-S$ operator, and $\pi$,
respectively (the '$-$' sign in front of $S$ is due to the fact that,
given a short exact sequence, and shifting by one, the boundary
induced in cohomology is '$-$ the initial boundary'; it also explains
the '$-$' sign in our definition of $\beta$).  Hence, by naturality,
$\pi= (-S)(-\beta)= S\beta$. Similar arguments
show that $p_2 r$ induces a chain map $\beta'$ between the kernels of
$p_3$ and $p_4$. By a diagram chasing and (\ref{sumunu}), we have
$d\compose \beta'\compose p_5= d\compose \beta\compose p_5$, hence
 $\beta'=\beta$. Using the naturality of the long exact
 sequences induced in cohomology by the right vertical and the bottom
 horizontal, we find $\pi= \beta S$.\\
Hence we are left with proving that $\beta$ induces isomorphism in
cohomology. We define now a new map on our diagram:
\[\tilde{\alpha}:= -\tilde{\phi}_1\tilde{d}: (I^{\ps n- 1\pd}/Im d, b_0)[1] \rmap
(I^{\ps n\pd}_{\natural} ,\del_0). \]
First of all, $-\tilde{\alpha}N= \tilde{\phi}_1\tilde{d}N=
\tilde{\phi}_1Nd= d$ by (\ref{uli}), so $\tilde{\alpha}$ induces a map
between the Cokernels of $N$ and $d$:
\[ \alpha: (I^{\ps n- 1\pd}/Im d+ Im N , b_0)[1] \rmap (I^{\ps
  n\pd}_{\natural}/Im d, \del_0) .\]
Since $\phi_1d\theta N= \phi_1(1- \theta d) N$, $\phi_1 N \equiv 1$ modulo $Im(1-
t)$, and $\phi_1 \theta= 0$ (cf. \ref{formulas}), we have:
\[ \alpha \compose \beta\, =\, 1\ .\]
Now we show that $\beta\alpha = 1$ in cohomology. Since $\theta N \phi_1 d= \theta (1+ \phi_0 (1- t)) d\equiv 1+
\theta \phi_0 (1- t) d$ modulo $Im(d)$, it suffices to show that:
\[ \theta \phi_0 d(1- t): (I^{\ps n- 1\pd}/Im d+ Im N , b_0)\rmap
(I^{\ps n- 1\pd}/Im d+ Im N , b_0) \]
is trivial in cohomology. For this we remark that our map factors as:
\[  (I^{\ps n- 1\pd}/Im d+ Im N , b_0)\stackrel{1- t}{\longrightarrow} 
(I^{\ps n- 1\pd}/Im d , \del_0)\stackrel{\theta \phi_0 d}{\longrightarrow}(I^{\ps
  n- 1\pd}/Im d+ Im N , b_0), \]
where the second map is a chain map by the non-trivial \ref{formulas}
  (iv), and the
  middle complex is contractible (by the usual $s_{-1}$).\\
Now, using Lemma \ref{imp}, and Lemma \ref{twist}, it is easy to see
  that all the formulas and arguments localize without any problem;
 this concludes the proof of  Theorem \ref{comp}. 

\begin{num}{\bf Example (Equivariant cycles):}\emph{ We point out that the new complexes
  computing $HC_{\delta}^*(\Ha)$, arising from Lemma \ref{ncom},
  appear naturally in the construction of characteristic maps
  associated to equivariant cycles. Recall \cite{Co3} that a chain of
  dimension $n$ is a triple $(\Omega, d, \int)$ where $\Omega=
  \oplus_{j=0}^{n} \Omega^n$ is a DG algebra, and $\int: \Omega^n\rmap
  {\bf C}$ is a graded trace on $\Omega$. It is a cycle if $\int$ is
  closed. If $(\Omega, d)$ is a $\Ha$ DG algebra, and $\int$ is
  $\delta$-invariant, we say that $(\Omega, d, \int)$ is a
  $ \Ha$-chain. For instance, if a Lie group $G$ acts smoothly on
  the $\Omega^j$'s, $d(g \omega)= g d(\omega)$, $\int
  g\omega= \int \omega$, for $g\in G$, $\omega\in \Omega$ (i. e.
  $(\Omega, d, \int)$ is a $G$- equivariant chain), then, with the
  induced infinitesimal action , $(\Omega, d,
  \int)$ is an $U(\lie)$-chain.\\
A $\Ha$-cycle over an algebra $A$ is given by such a cycle, together
  with an algebra homomorphism $\rho: A\rmap \Omega^0$. We view $\rho$
  as an element of (total) degree $1$ on the bigraded differential
  algebra $L= Hom(BA, \Omega)$. The structure on $L$ is the one
  induced by the graded structures on $B(A)$, and $\Omega$,
  respectively. That is, the bigrading: $L^{p, q}= Hom(A^{\otimes p},
  \Omega^q)$, the differentials $d(f)= d\compose f$, $\del_0(f)=
  -(-1)^{deg\ps f\pd}f\compose b\,'$, the product: $(\phi *\psi)(a_1,
  .\, .\, .\, , a_{p+q})=$ 
$(-1)^{p deg\ps\psi\pd}
  \phi(a_1, .\, .\, .\, , a_{p})$ $\psi(a_{p+1}, .\, .\, .\, , a_{p+q})$.
As in \ref{caract}, there is a unique $\Ha$ DG algebra map $k:
  W(\Ha)\rmap Hom(BA, \Omega)$, sending $1$ to $\rho$; it is
  compatible with the bigraded differential structure. On the other
  hand, using $\int: \Omega^n\rmap {\bf C}[n]$, we have a graded trace
  $\int^{\natural}: Hom(BA, \Omega)\rmap C_{\lambda}(A)[n]$, similar
  to (\ref{buc2}), $\int^{\natural}(f)= \int\compose f\compose N$,
  whose composition with $k$ is still denoted by $k: W(\Ha)\rmap C_{\lambda}(A)[n]$.\\
Now, since $\int$ is closed, $\delta$-invariant, and supported on
degree $n$, the relevant target of $k$ is the complex
$(I_{\natural, \delta}^{\ps n\pd}/Im d, \del_0)$, hence it induces a
characteristic map
$k: HC_{\delta}^{*- 2n- 1}(\Ha)$ $\cong H^*(I_{\natural, \delta}^{\ps
  n\pd}/Im d, \del_0)$ $\rmap H^*(C_{\lambda}(A)[n+1])$ $\cong HC^{*- n-
  1}(A$). As in the proof of \ref{damm}, it is compatible with the $S$-operator.\\
Let's assume now that $(\Omega, d, \int)$ is cobordant to the trivial
cycle, that is, there is an $n$-dimensional chain $(\tilde{\Omega}, \tilde{d},
\tilde{\int})$, $\tilde{\rho}: A\rmap \tilde{\Omega}^0$, and an
equivariant chain map $r: \tilde{\Omega}\rmap\Omega$ such that
$\tilde{\int}\compose \tilde{d}= \int\compose r$, $r\compose
\tilde{\rho}=\rho$. As before, we have an induced map $\tilde{k}:
(I_{\natural, \delta}^{\ps n+1\pd}/Im d, \del_0)[-1]\rmap
C_{\lambda}(A)[n+1]$. Since $\tilde{\int}\compose \tilde{d}=
\int\compose r$, we have $\tilde{k}\compose d= k$, and we are on the
localized version of the short exact sequence (\ref{nustiu}). Since
its boundary identifies with $\pi$, hence with the $S$ operator cf. \ref{comp}, we
deduce that, after stabilizing by $S$, $k$ is trivial (in
  cohomology). We summarize our discussion:}\end{num}

\begin{cor}\label{las} Let $A$ be an algebra. The previous construction
  associates to any $n$-dimensional equivariant cycle $(\Omega, d,
  \int)$ over $A$ a characteristic map compatible with the
  $S$-operator:
\[ k: HC_{\delta}^{*}(\Ha)\rmap HC^{*+ n}(A)\ .\]
If $k_0, k_1$ are associated to cobordant cycles, then $S\compose k_0=
S\compose k_1$.
\end{cor}

For instance, for $\Ha= U(\lie)$ previously mentioned, using Theorem
\ref{env} we get:

\begin{cor} Given a Lie group $G$, and  a smooth $n$-dimensional
  $G$-cycle over an algebra $A$, we have induced  ($\bf{Z}/2\bf{Z}$ graded)
  characteristic maps:
\[ k: H_{*}(\lie)\rmap HP^{*+ n}(A) \ .\]
Cobordant $G$-cycles induce the same map.
\end{cor}

\end{document}